\magnification=1100
\overfullrule=0mm

\input xy \xyoption{all}

\font\fiverom=cmr10 at 5pt
\font\sevenrom=cmr10 at 7pt
\font\eightrom=cmr10 at 8pt
\font\ninerom=cmr10 at 9pt

\font\twelvecmb=cmb10 at 12pt
\font\tencmb=cmb10 at 10pt
\font\ninecmb=cmb10 at 9pt

\font\sevencmb=cmb10 at 7pt

\font\nineromtt=cmtt10 at 9pt

\font\bb= msbm10
\font\sevenbb= msbm10 at 7pt





\font\eightronde=cmsy10 at 8pt

\def\HC{\hbox{$\cal H$}}
\def\HCp{\hbox{\eightronde H}}

\def\HBB{\hbox{\bb H}}
\def\HBBp{\hbox{\sevenbb H}}
\def\Hr{\hbox{\rm H}}

\def\Hom{\hbox{$\rm Hom$}}

\def\Aut{{\rm Aut}}
\def\End{\mathop{\rm End}\nolimits}
\def\Map{\mathop{\rm Map}\nolimits}
\def\Tors{{\rm {\cal T}\! ors}}
\def\Torsr{{\rm Tors}}
\def\tors{{\rm tors}}

\def\relbar{\mathrel{\smash-}}

\def\lodot{ \mathop{\circ\mkern-7.05mu\cdot\mkern3.5mu}\nolimits }

\def\ZBIc{\hbox{\sevencmb Z}}

\def\Zb{\hbox{\tencmb Z}}
\def\ZBB{\hbox{\bb Z}}

\def\ZC{\hbox{$\cal Z$}}

\def\AC{\hbox{$\cal A$}}
\def\BC{\hbox{$\cal B$}}
\def\CC{\hbox{$\cal C$}}

\def\T{\hbox{$\cal T$}}

\def\Cr{\hbox{\rm C}}
\def\Er{\hbox{\rm E}}
\def\Hr{\hbox{\rm H}}
\def\Wr{\hbox{\rm W}}
\def\Zr{\hbox{\rm Z}}

\def\leqp{{\scriptscriptstyle\leq}}
\def\Del {\hbox{$\Delta$}}
\def\id{\hbox{\rm id}}
\def\idr{\hbox{\sevenrom id}}
\def\eps {\hbox{$\varepsilon$}}

\def\qed{\lower 2pt\hbox{\vrule\vbox to 5.5pt{\hrule width 4.5pt
                                             \vfill
                                             \hrule}\vrule}}

\def\dm{\unskip\hfill\hfill\penalty 500\kern 10pt\qed}

\def\lr {\hbox{$\ \longrightarrow\ $}}
\def\ot {\hbox{$\otimes$}}
\def\pa{\S\kern.15em}
\def\Dem{\noindent {\sl Proof:}$\, \,$}

\def\longrighthook{\lhook\joinrel\relbar\joinrel\rightarrow}

\def\hfll#1#2{\smash{\mathop{\hbox to 10mm{\rightarrowfill}}
\limits^{\scriptstyle#1}_{\scriptstyle#2}}}

\def\cross{{\times}}
\def\lcross{ {\times}}
\def\llcross{{\times}}

\newtoks\auteurcourant      \auteurcourant={\hfil}
\newtoks\titrecourant       \titrecourant={\hfil}

\newtoks\hautpagetitre      \hautpagetitre={\hfil}
\newtoks\baspagetitre       \baspagetitre={\hfil}

\newtoks\hautpagegauche   \newtoks\hautpagedroite 
  
\hautpagegauche={\eightrom\rlap{\folio}\eightrom\hfil\the\auteurcourant\hfil}
\hautpagedroite={\eightrom\hfil\the\titrecourant\hfil\eightrom\llap{\folio}}

\newtoks\baspagegauche      \baspagegauche={\hfil} 
\newtoks\baspagedroite      \baspagedroite={\hfil}

\newif\ifpagetitre          \pagetitretrue

\headline={\ifpagetitre\the\hautpagetitre
            \else\ifodd\pageno\the\hautpagedroite
             \else\the\hautpagegauche
              \fi\fi}

\footline={\ifpagetitre\the\baspagetitre\else
            \ifodd\pageno\the\baspagedroite
             \else\the\baspagegauche
              \fi\fi
               \global\pagetitrefalse}

\def\raggedbottom{\topskip 10pt plus 36pt\r@ggedbottomtrue}

\newcount\notenumber \notenumber=1
\def\note#1{\footnote{$^{{\the\notenumber}}$}{\eightrom {#1}}%
\global\advance\notenumber by 1}

\auteurcourant={Marc WAMBST}
\auteurcourant={Philippe NUSS -- Marc WAMBST}

\titrecourant={NON-ABELIAN HOPF COHOMOLOGY}

\null
\vfill

\noindent {\twelvecmb NON-ABELIAN HOPF COHOMOLOGY II

}
\smallskip
\smallskip
\noindent {\tencmb -- THE GENERAL CASE  --}
\vskip 10pt

\bigskip

\noindent {PHILIPPE NUSS, MARC WAMBST}

\vskip 3pt
\noindent {\ninerom Institut de Recherche Math\'ematique Avanc\'ee,
Universit\'e Louis-Pasteur et CNRS, 7, rue Ren\'e-Descartes,
67084 Strasbourg Cedex, France. 
e-mail: {\nineromtt nuss@math.u-strasbg.fr} and 
{\nineromtt wambst@math.u-strasbg.fr}}

\vskip 25pt
\itemitem{}{\ninecmb Abstract.} {\ninerom 
 We introduce and study non-abelian cohomology sets of Hopf algebras with coefficients in Hopf comodule algebras.
We prove that these sets generalize as well Serre's non-abelian group cohomology theory as
the cohomological theory constructed by the authors in a previous article. We esta\-blish their functoriality and compute explicit examples. Further we classify Hopf torsors.}

\vskip 5pt
\itemitem{}\noindent {\ninecmb MSC 2000 Subject Classifications.} 
{\ninerom Primary: 18G50, 16W30,  14A22; Secondary:   20J06, 55U10.} 

\vskip 5pt
\itemitem{}\noindent {\ninecmb Key-words:} {\ninerom non-abelian cohomology, 
Hopf comodule algebra, torsor,  cosimplicial non-abelian groups.}

\vskip 30pt  

\noindent {\tencmb I{\ninecmb NTRODUCTION}.}   The present article, conceived as the continuation of [8], is devoted to the study of non-abelian cohomology theory in the  Hopf algebra setting. 
We define a general  cohomology theory analogous to that for groups  adapted to Hopf algebras and suitable coefficient objects  ([9], [10]).
In~order to clarify the purpose of our work, we recall first some basic facts about the classical constructions in the framework if groups.
Let $G$ be a  group  acting on a group~$A$. The non-abelian cohomology theory $\Hr^*(G, A) $ for groups  may be organized in three different stages depending on the properties of the coefficient group~$A$.

\medskip
 \itemitem{Stage 1)} {\sl  The group of coefficients $A$ is abelian.} The   classical  Eilenberg-MacLane
cohomology theory  $\Hr^*(G, A) = {\rm Ext}^*_{{\hbox {\ZBIc}} [G]}(\Zb, A)$ produces  a sequence of commutative groups. It provides useful invariants
in  homological algebra, algebraic topology and algebraic number theory. 
\smallskip
\itemitem{Stage 2)} {\sl  The group of coefficients $A$ is not abelian.}  The previous construction  fails in this case. However it is still possible  to define a group $\Hr^0(G, A)$ and a pointed set $\Hr^1(G, A)$. 
This theory, called the non-abelian cohomology theory of groups,
was introduced by Lang and Tate ([4])
for Galois groups with coefficients in an algebraic group, and
was studied in full generality by Serre ([9], [10]). It is for instance well-known that the non-abelian cohomology set $\Hr^1(G, A)$ classifies
the $G$-torsors (or principal homogeneous spaces)  on $A$ (see [10]).

\smallskip

\itemitem{Stage 3)} {\sl The group of coefficients $A$ is the group of automorphisms of a $G$-Galois extension.} 
Suppose that the group $G$ is  finite and acts as a Galois group on
a Galois extension $S/R$ of  noncommutative rings (for this generalization of Galois extensions of fields, see  [5]).
Let $M$ be an $S$-module endowed with a compatible $G$-action. The latter induces a $G$-group structure on the group $\Aut _S(M)$ of $S$-linear automorphisms of $M$.
One of the authors ([7]) showed  that in this context  non-abelian cohomology theory comes into play. In particular he proved that the set $\Hr ^1(G, \Aut _S(M))$ classifies objects which arise in descent theory
along $S/R$, for example descent cocycles on $M$
or twisted forms of~$M$.

\medskip 
\eject

Hopf algebras naturally generalize groups. Kreimer and Takeuchi ([3]) widened  
Galois extensions to  Hopf-Galois extensions of rings in the following spirit.
As a group acting on rings  plays the r\^ole  of the symmetry object for Galois extensions,  a Hopf algebra coacting on the rings does for Hopf-Galois extensions.
In [8], we answered the 
natural question of extending Stage 3 to this setting. 
For a Hopf algebra $H$,  an $H$-Hopf comodule algebra $S$, and  an $(H,S)$-Hopf module $M$,    we introduced a group  $\Hr ^0 (H, M)$
and a pointed set  
$\Hr ^1 (H, M)$.  This construction, here called {\sl restricted non-abelian Hopf cohomology theory}, replaces $\Hr ^*(G, \Aut _S(M))$. It offers a generalization of Stage 3 in the following two senses
(see [8]):

\item{--} If $S/R$ is an $H$-Hopf-Galois extension, then $\Hr ^1 (H, M)$ classifies the analogue of descent cocycles on $M$ along $S/R$ and  the twisted forms of $M$. 

\item{--} Given a group  $G$, a $G$-Galois extension of rings $S/R$ is nothing but a $\Zb^G$-Hopf-Galois
extension, where $\Zb^G$ stands for the Hopf algebra of functions on $G$. Then $\Hr ^* (\Zb^G, M)$ is isomorphic to~$\Hr ^*(G, \Aut _S(M))$.

\medskip
The aim of this article is to define a non-abelian cohomology theory 
in the Hopf context corresponding to Stage 2. More precisely, let $H$ be a Hopf algebra over a commutative ring $k$. For any $H$-comodule $k$-algebra $E$ we
introduce the {\sl general non-abelian Hopf cohomology theory of $H$ with coefficients in $E$}.
We define  a group $\HC^0(H, E)$ and a pointed set $\HC^1(H, E)$. Theses constructions are based on the non-abelian cohomology  theory associated to a
pre-cosimplicial group.
We prove three main results (the precise wording  and definitions will be found in the core of the article):

\item{(a)} We show  (Theorem 1.5) that the cohomology theory $\HC^*(k^G, E)$
is isomorphic to $\Hr ^*(G, E^{\lcross})$, where~$k^G$ denotes the Hopf algebra of functions on $G$ and $E^{\lcross}$ is the group of invertible elements of $E$. 

\item{(b)} Let
$S$ be an $H$-Hopf comodule algebra and $M$ be an $(H,S)$-Hopf module. We establish (Theorem 2.6) that under lax technical conditions,
 $\HC^*(H, \End _S(M))$ and $\Hr^*(H, M)$ are isomorphic.

\item{(c)} Finally, if $E$ is an $H$-comodule algebra, we classify  $(H,E)$-Hopf torsors via the pointed set~$\HC^1(H, E)$ (Theorem 3.4).

\medskip

The article is built in the following way. The first section is devoted to the definition and the properties of general non-abelian Hopf cohomology theory. 
There we prove Result (a), give some examples, explicit computations   (\pa 1.2 and \pa 1.4), and show that the 
Hopf module structures may be deformed with the help of  $1$-cocycles (Proposition 1.7). In \pa 1.6 we study the functoriality of the general non-abelian Hopf cohomology sets
and write down an exact sequence associated to a sub-comodule algebra.
In the second section, we clarify the links between general and restricted non-abelian Hopf cohomo\-logy theory. To this end, we state a technical condition (Condition $({\cal F}_n)$
in \pa 2.2) which allows to endow 
the endomorphism algebra of an Hopf module with a comodule structure (Lemma 2.4). We then deduce Result (b).
The third and last section deals with Hopf torsors. We define them as a generalization of usual torsors  (Definition 3.2, Proposition 3.7, and Corollary 3.8) and prove Result (c).

\medskip
We mention here that an attempt of generalizing the non-abelian group cohomology theory to the Hopf context was done by Blanco Ferro
([1]). This author adapted Sweedler's theory 
([11]), which can be viewed as a generalization of Stage 1. Blanco Ferro defined a $1$-cohomology set $\Hr^1(H, A)$, where $H$ is a cocommutative Hopf algebra and $A$
is an algebra not necessarily commutative. 
His construction is in some sense dual to ours. But if one tries to apply it to the Hopf-Galois extensions, one has to restrict
oneself to a very particular case: 
not only does $H$ have to be a commutative finitely generated $k$-projective Hopf algebra,
but the Hopf-Galois extension  $S/k$ is over the gound field and moreover has to be  commutative.

\eject

\medskip
\medskip
\noindent {\tencmb 0. Conventions, notations, and terminology.}

\smallskip

\noindent 
Let $k$ be a fixed commutative and unital ring. 
The unadorned symbol $\ot$ between a right {\hbox{$k$-module}} and a left
$k$-module stands for $\ot _k$. 
By {\sl (co-)algebra} we mean a (co-)unital (co-)associative {\hbox{$k$-(co-)}}al\-gebra.
By {\sl (co-)module} over a (co-)algebra $D$, we always understand a right $D$-(co-)module unless otherwise stated.
Let $M$ be a $k$-module. We identify in a systematic way $M \ot k$ with $M$.

For any algebra $D$, we denote by $D^{\lcross}$ the group of invertible elements in $D$.
If $M$ is a $D$-module, $\End _D(M)$ (respectively $\Aut _D(M)$) is the algebra (respectively the group) of $D$-linear
endomorphisms (respectively automorphisms) of $M$.

Let $H$ be a Hopf algebra  with multiplication $\mu_H$, unity map $\eta _H$,
comultiplication $\Delta _H$, counity map $\varepsilon _H$, and  antipode $\sigma_H$.
Recall that an $H$-comodule algebra $E$ is a $k$-module which is both  an algebra and an $H$-comodule such that
the coaction map  is a morphism of algebras. A morphism of $H$-comodule algebras
is simultaneously a morphism of algebras and of $H$-comodules.
Suppose that $E$ is an $H$-comodule algebra.
Let $M$ be both an $E$-module 
and an $H$-comodule. If the coaction map
$\Delta _M: M \lr M \ot H$  verifies the equality
$$\Delta _M (ms) = \Delta _M(m) \Delta _S(s)$$ for any $m \in M$ and $s \in E$, we say that
$M$ is an {\sl $(H,E)$-Hopf module} 
(also called a {\sl relative Hopf module} in the literature)  and that
$\Delta _M$ is  {\sl $(H,E)$-linear}.
A {\sl morphism of $(H,E)$-Hopf modules} is an $E$-linear map
 $f: M \lr M'$ such that $(f \ot \id _M) \circ\Delta _M = \Delta _{M'} \circ f$.
Observe that  $E$ itself is naturally an $(H,E)$-Hopf module.

To denote the coactions on elements, we use the Sweedler-Heyneman convention, that is, for $m \in M$, we write
 $\Delta _M(m) = m_0 \ot m_1$, with summation implicitly
understood. More generally, when we write down a tensor we usually omit the summation sign $\sum$.

\medskip

Let $G$ be a finite group with neutral element $e$. Denote by $k^G$ the $k$-free Hopf algebra over the $k$-basis
$\{ \delta _g\} _{g \in G}$, with the following structure maps:
the multiplication is given by ${\displaystyle \delta _g \cdot \delta _{g'} = \partial _{g,g'} \delta _{g}}$,
where $\partial _{g,g'}$ stands for the 
Kronecker symbol of $g$ and $g'$; 
the  comultiplication $\Del _{k^G}$ is
defined by 
${\displaystyle \Del _{k^G} (\delta _g) =
\sum _{ab = g}\delta _a \ot \delta _b}$;
the unit in $k^G$ is the element 
${\displaystyle 1 = \sum _{g \in G}\delta _g}$;
the counit $\varepsilon _{k^G}$ is
defined by $\varepsilon _{k^G} (\delta _g) = \partial _{g,e}1$; 
the antipode $\sigma_{k^G}$ sends $\delta _g$ on $\delta _{g^{-1}}$.
When $k$ is a field, then $k^G$ is the dual of the usual group algebra $k[G]$.

\bigskip
\smallskip\smallskip

\noindent {\tencmb 1. General non-abelian Hopf cohomology theory.}

\smallskip
\noindent The first section is devoted to the definition, the properties and examples of general non-abelian Hopf cohomology theory.
The constructions are provided in simplicial terms (a r\'esum\'e 
about the simplicial language may be found in [6]).

\medskip

\noindent {\sl 1.1. Definitions.} Let $\AC ^* =  \xymatrix{ A^0  \ar@<1.3ex>[r]^{d^0} 
\ar@<-1.3ex>[r]^{d^1}  & \ A^1
\ar@<2.3ex>[r]^{d^0} 
\ar@<0ex>[r]^{d^1} 
\ar@<-2.3ex>[r]^{d^2}& \ A^2}$ be a  pre-cosimplicial group. The {\sl non-abelian $0$-cohomology 
group} $\HBB^0(\AC ^*)$ is the equalizer
of the pair $(d^0, d^1)$:
$$\HBB^0(\AC ^*) = \{ x \in A^0 \ \vert \ d^1(x) = d^0(x) \}.$$ 
 The {\sl non-abelian $1$-cohomology pointed set} $\HBB^1(\AC ^*)$ is the right quotient
$$\HBB^1(\AC ^*) =  A^0 \backslash \ZBB^1(\AC ^*).$$ 
Here the set $\ZBB^1(\AC ^*)$ of {\sl $1$-cocycles} 
is the subset of $A^1$ defined by 
$$\ZBB^1(\AC ^*) =\{ X \in A^1 \ \vert \ d^2(X)d^0(X) = d^1(X) \}.$$ 
The group $A^0$
acts on  the right on $A^1$  by $$X \leftharpoonup x = (d^1x^{-1}) X (d^0x),$$
where $X \in A^1$ and $x \in  A^0$. Using the pre-cosimplicial relations, one easily checks that this action restricts to  $\ZBB^1(\AC ^*)$.
Two $1$-cocycles $X$ and $X'$
are said to be {\sl cohomologous} if they belong to the same orbit
under this action. The quotient set $\HBB^1(\AC ^*) =   A^0 \backslash \ZBB^1(\AC ^*) $
is pointed with distinguished point the class of the neutral element of $  A^1$.
\medskip

\goodbreak
Let $H$ be a Hopf algebra,
let $E$ be an $H$-comodule algebra with multiplication $\mu_E$ and coaction
$\Delta _E$. 
We define two maps $d^i: E \lr E \ot H$
($i= 0,1$)
and three maps $d^i:  E \ot H \lr  E \ot H\ot H$ ($i=0,1,2$) 
by the formulae
$$\eqalign{d^0 (x)  &  = \Delta _E (x), \hskip47.55pt d^1 (x)  = x\ot 1,  \cr
 d^0(X)  & = (\Delta _E \ot \id _H) (X), \quad
 d^1(X)    = (\id _E \ot  \Delta _H)(X),  \quad   \ d^2(X)    = X \ot 1, \cr}$$
where $x \in E$ and $X \in E \ot H$.
\goodbreak
\medskip
\noindent {\tencmb Lemma 1.1.} {\sl The  diagram ${\cal C}_{\leqp2}(H, E)$ given by
$$  \xymatrix{ E  \ar@<1.3ex>[r]^{d^0 \ \ \ } 
\ar@<-1.3ex>[r]^{d^1 \ \ \ }  & \ E \ot H
\ar@<2.3ex>[r]^{d^0\ \ \  } 
\ar@<0ex>[r]^{d^1\ \ \ } 
\ar@<-2.3ex>[r]^{d^2\ \ \ }& \ E \ot H \ot H}$$ 
is a pre-cosimplicial object in the category of algebras.}

\medskip

\Dem The maps $d^i$ are easily seen to be morphisms of algebras. 
The pre-cosimplicial relations $d^id^j = d^jd^{i-1}$ for $i > j$ follow from the  Hopf axioms for 
$H$  and $E$.
 \dm

\medskip

Lemma 1.1 allows us to deduce a  pre-cosimplicial diagram ${\cal C}^{\lcross}_{\leqp2}(H, E) $ in the category of groups  by setting:

$$ \xymatrix{ E^{\lcross}  \ar@<1.3ex>[r]^{d^0 \ \ \ } 
\ar@<-1.3ex>[r]^{d^1 \ \ \ }  &  \ (E \ot H)^{\lcross}
\ar@<2.3ex>[r]^{d^0 \ \ \ } 
\ar@<0ex>[r]^{d^1 \ \ \ } 
\ar@<-2.3ex>[r]^{d^2\ \ \ }& \ (E \ot H \ot H)^{\lcross}}$$
(we still denote by  $d^i$ the restrictions  of the maps $d^i : E \ot H^{\otimes j} \lr E \ot H^{\otimes (j+1)}$ to the corresponding multiplicative groups).

\medskip

\noindent {\sl Remark~:} Both ${\cal C}_{\leqp2}(H, E) $ and ${\cal C}^{\lcross}_{\leqp2}(H, E) $ are in fact cosimplicial objects. The 
codegeneracy maps on ${\cal C}_{\leqp2}(H, E) $  are given by
$$\eqalign{s^0 & = \id _E \ot \varepsilon _H : E\ot H \lr E, \cr s^0 & = \id _E \ot \varepsilon _H \ot \id _H : E\ot H \ot H \lr E \ot H \quad {\hbox{\rm and}} \quad 
s^1 = \id _E \ot \id _H \ot \varepsilon _H : E\ot H \ot H \lr E \ot H.\cr}$$
The codegeneracy maps on ${\cal C}^{\lcross}_{\leqp2}(H, E) $ are again obtained by restriction.

\medskip

\goodbreak

\noindent {\sl Definition 1.2~:} The {\sl general non-abelian Hopf cohomology objects} $\HC^*(H, E) $ of a Hopf algebra $H$ with coefficients in an $H$-comodule algebra $E$
is  the non-abelian cohomology theory associated to the pre-cosimplicial diagram ${\cal C}^{\lcross}_{\leqp2}(H, E) $.

\medskip
In other words 
$$\eqalign{\HC^0(H, E) & = \HBB^0 \bigl({\cal C}^{\lcross}_{\leqp2}(H, E) \bigr) = \{ x \in E^{\lcross} \ \vert \ d^1(x) = d^0(x) \} \quad {\hbox{\rm  and}}\cr
\HC^1(H, E) &= \HBB^1 \bigl({\cal C}^{\lcross}_{\leqp2}(H, E) \bigr) = E^{\lcross} \backslash \ZC^1(H, E).\cr}$$
Observe that $\HC^0(H, E)$ is the group $(E^{{\eightrom co}H})^{\lcross}$ of invertible coinvariant elements of $E$.
The set $\ZC^1(H, E)$ of {\sl Hopf $1$-cocycles of $H$ with coefficients in $E$} 
is the subset of $(E\ot H)^{\lcross}$ given by 
$$\ZC^1(H, E) =\{ X \in (E\ot H)^{\lcross} \ \vert \ d^2(X)d^0(X) = d^1(X) \}.$$ 
We refer to  $d^2(X)d^0(X) = d^1(X)$ as the {\sl Hopf $1$-cocycle relation}. 

\goodbreak
\medskip

\noindent {\sl Remarks 1.3~:} 

\noindent a) For any Hopf algebra $H$ and any $H$-comodule algebra $E$, one proves the inclusion $$\ZC^1(H, E) \subseteq 
{\rm Ker}\bigl(\id _E\ot \varepsilon _H: (E \ot H)^{\lcross} \lr E^{\lcross}\bigr)$$ by applying the map $\id _E \ot \varepsilon _H \ot \id _H$ to 
the  Hopf $1$-cocycle relation.

\noindent b) If the algebras $E$ and $H$ are both commutative, the sets $\ZC^1(H, E)$ and $\HC^1(H, E)$ become groups with product induced by the multiplication of $E \ot H$.

\medskip

\goodbreak
\noindent {\sl 1.2. First examples.} 

\noindent 1) {\sl The Hopf algebra is trivial.} Any algebra $E$ is naturally a $k$-comodule algebra with the coaction $\Delta _E$ equal to $\id _E$. One then has:
$$\HC^0(k, E) = E^{\lcross}  \quad {\hbox{\rm  and}} \quad \HC^1(k, E) = \{1\}.$$ Indeed, the first equality is obvious. One checks that
$\ZC^1(k, E)$ is the pointed set of invertible idempotent elements of $E$, that is nothing else than $\{1\}$.

\medskip

\noindent 2) {\sl The coefficients are trivial.} Let $H$ be a Hopf algebra. The  ground ring $k$ is an $H$-comodule algebra through the coaction  $\Delta _k$ given by the unity map $ \eta _H$. Denote by ${\rm Gr}(H)$ the group 
of grouplike elements in $H$.
One then has: $$\HC^0(H, k) = k^{\lcross}  \quad {\hbox{\rm  and}} \quad \HC^1(H, k) \cong {\rm Gr}(H),$$ the latter relation being an isomorphism of groups.
The calculation of $\HC^0(H, k)$ is straightforward. We compute now $\ZC^1(H, k)$. A $1$-cocycle is in particular an element
$h\in H$ verifying the $1$-cocycle relation, here  $h \ot h = \Delta _H(h)$. So the element $h$ is grouplike, hence incidentaly also invertible in $H$.
The action of $k^{\lcross} $ on $\ZC^1(H, k)$ is trivial; therefore $\HC^1(H, k)$
is the whole group of grouplike elements of $H$.

\medskip

\noindent 3) {\sl The coefficients are the Hopf algebra itself.}  A Hopf algebra $H$ is a comodule algebra over itself. One has:
$$\HC^0(H, H)  = k^{\lcross} \quad {\hbox{\rm  and}} \quad   \HC^1(H, H) = \{1\}.$$
The first equality follows from  the very definition: $\HC^0(H, H) = (H^{{\eightrom co}H})^{\lcross}$. To prove the second equality,
pick $X \in \ZC^1(H, H)$ and apply the map
$\varepsilon _H \ot \id _H \ot \id _H$ to the cocycle relation $d^2(X )d^0(X) = d^1(X)$. One gets
$(x \ot 1)X = \Delta _H(x)$, with $x = (\varepsilon _H \ot \id _H) (X)$. So $\ZC^1(H, H)$ is contained in the set $\{(x^{-1} \ot 1)\Delta _H(x) \ \vert \ x \in H^{\lcross}\}$, which is equal to  
$\{d^1(x^{-1})d^0(x) \ \vert \ x \in H^{\lcross}\}$.
Conversely, if $X = (x^{-1} \ot 1)\Delta _H(x)$ for $x \in  H^{\lcross}$, then $X$ fulfills the cocycle relation. So $\ZC^1(H, H)$ equals  
$\{d^1(x^{-1})d^0(x) \ \vert \ x \in H^{\lcross}\}$, and therefore the $1$-cohomology set is trivial.

\goodbreak
\bigskip
\noindent {\sl 1.3. Link with non-abelian group cohomology.} 
We first recall the definitions given by Serre ([9], [10]) of the non-abelian cohomology theory $\Hr^i(G,A) $ (with $i = 0,1$) of a
 group $G$ with coefficients in  a
(left) $G$-group $A$.
The 0-cohomology object
$\Hr^0(G,A) $ is the group $A^G$ of 
invariant elements of $A$ under the action of $G$.
The set $\Zr^1(G, A)$ of $1$-cocycles is given by
$$\Zr^1(G, A) = \{ \alpha : G \lr A \ \vert \ \ \alpha (gg') =
\alpha (g){{}^{g}\! \bigl(}\alpha (g')\bigr),
 \ \ \forall \ g, g' \in G \}.$$
It is pointed with distinguished point the constant map $1: G \lr A$.
The group
$A$ acts on the right on $\Zr^1(G, A)$  by $$(\alpha \leftharpoonup a)(g) = a^{-1}\alpha (g) \ {}^{g}\! a,$$
where $a \in A$, $\alpha  \in \Zr^1(G, A)$, and $g \in G$.
Two $1$-cocycles $\alpha $ and $\alpha '$
are {\sl cohomologous} if they belong to the same orbit under this action.
The non-abelian $1$-cohomology set  $\Hr^1(G, A)$ is the left quotient
$A \backslash \Zr^1(G, A)$.
It is pointed with distinguished point the class of the constant map $1: G \lr A$.
\medskip

The non-abelian cohomology theory of groups may be interpreted as the non-abelian cohomology theory associated to the pre-cosimplicial diagram of groups
$${\cal G}_{\leqp2}(G, A)  = \Bigl(\xymatrix{ A = \Map (G^0,A) \ \ar@<1.3ex>[r]^{\ \ \  d^0 }
\ar@<-1.3ex>[r]^{\ \ \  d^1  }  & \ \Map (G,A)
\ar@<2.3ex>[r]^{d^0\ \ \  } 
\ar@<0ex>[r]^{d^1\ \ \ } 
\ar@<-2.3ex>[r]^{d^2\ \ \ }& \ \Map (G^2,A)\Bigr).}$$ 
Here $\Map(G^i,A)$ stands for the set of the maps from $G^i$ to $A$, which is  endowed with the group structure induced by pointwise multiplication.
The coboundaries are given by $$\eqalign{d^0(x) &: g \longmapsto {^g\!x}  , \hskip47.55pt d^1(x) : g \longmapsto  x,  \cr
 d^0(\alpha )  & : (g,g') \longmapsto {^g\!}\alpha (g'), \quad
 d^1(\alpha ) : (g,g') \longmapsto \alpha (gg'),  \quad   \ d^2(\alpha ) : (g,g') \longmapsto \alpha (g), \cr}$$
where $x \in A$, $g,g'\in G$ and $\alpha \in \Map (G,A)$. The reader may easily check that the pre-cosimplicial relations are satisfyied and that one has the equality
$$\HBB^*({\cal G}_{\leqp2}(G, A) ) = \Hr^*(G,A).$$

\medskip

We now connect the general non-abelian Hopf cohomology theory with the  non-abelian cohomology theory of groups.
Let $G$ be a finite group and $E$ be a  $k^G$-comodule algebra. For any $x\in E$, write
$$\Delta_E(x) =  \sum _{g \in G} {^g\!x} \ot \delta _g.$$
This formula defines an action of the  group $G$ on the algebra $E$, hence on the group $E^{\lcross}$.
One has the following result~:

\medskip

\noindent {\tencmb Proposition 1.4.} {\sl  Let $G$ be a finite group and  $E$ be a  $k^G$-comodule algebra. The   
pre-cosimplicial groups ${\cal G}_{\leqp2}(G, E^{\lcross})$ and ${\cal C}_{\leqp2}(k^G, E)$ are isomorphic.
}

\medskip
\noindent Before we give the proof, we state the following immediate consequence: 
\medskip

\goodbreak
\noindent {\tencmb Theorem 1.5.} {\sl  Let $G$ be a finite group and  $E$ be a  $k^G$-comodule algebra. There is the   equality of groups
$$\HC^0(k^G, E) = \Hr^0(G, E^{\lcross})$$ and an isomorphism of pointed sets
$$ \HC^1(k^G, E) \cong \Hr^1(G, E^{\lcross}).$$ }

\noindent {\sl Proof of Proposition 1.4.} First remark that any element in $(E\ot k^G)^{\lcross}$  is 
of the form $\displaystyle \sum_{g\in G}x_g\ot \delta _g$, where for all $g\in G$, the element $x_g$ belongs to $E^{\lcross}$.
In the same way any element in
 $(E\ot k^G\ot k^G)^{\lcross}$ is 
of the form $\displaystyle \sum_{g,g'\in G}x_{g,g'}\ot \delta _g\ot \delta _{g'}$,  where for all $g, g'\in G$, the element $x_{g,g'}$ belongs to $E^{\lcross}$.
We consider the map $\gamma_*: {\cal C}_{\leqp2}(k^G, E) \lr {\cal G}_{\leqp2}(G, E^{\lcross})$, given by
$$\eqalign {&\gamma _0 = \id _{E^{\llcross}} \cr
&\gamma _1(\sum_{g\in G}x_g\ot \delta _g) : u \longmapsto x_{u} \cr
&\gamma _2(\sum_{g,g'\in G}x_{g,g'}\ot \delta _g\ot \delta _{g'}) : (u,v) \longmapsto x_{u,v}, \cr}$$
for any $u,v \in G$. On each level, $\gamma _*$ is an isomorphism of groups since $(E \ot k^G)^{\lcross}$ (respectively $(E \ot k^G\ot k^G)^{\lcross}$) is isomorphic to  $(E^{\lcross})^{|G|}$
(respectively to  $(E^{\lcross})^{|G|^2}$).

\goodbreak
It remains to check that $\gamma _*$ is a morphism of pre-cosimplicial objects, in other words $\gamma _*$  verifies $\gamma_jd^i = d^i\gamma_{j-1}$
for any $1\leq j \leq 2$ and $0\leq i \leq j$. This is done by  direct computations. For example, set 
$$\nu = \gamma _2d^0(\sum_{g\in G}x_g\ot \delta _g) =  \gamma _2(\sum_{g,g'\in G} {^{g'}\!x_g}\ot \delta _{g'} \ot \delta _g). $$ So $\nu (u, v) = {^u\!}x_v$, for any $u,v \in g$.
Hence $\displaystyle \nu = d^0\gamma_1(\sum_{g\in G}x_g\ot \delta _g)$.
As an other  example, set 
$$\nu' = \gamma _2d^1(\sum_{g\in G}x_g\ot \delta _g) =  \gamma _2(\sum_{h,h'\in G} {x_{hh'}}\ot  \delta _{h} \ot \delta _{h'}). $$ So $\nu ' (u, v) = x_{uv}$, for any $u,v \in g$.
Hence $\displaystyle \nu' = d^1\gamma_1(\sum_{g\in G}x_g\ot \delta _g)$. We leave to the reader the three remaining computations.
\dm

\medskip
\medskip

\noindent {\sl Two direct applications of Theorem 1.5.} 
\smallskip
\noindent 1)  Let $G$ be a finite group. One may recover the isomorphism between the group ${\rm Gr}(k^G)$  of grouplike elements of $k^G$
and the Pontryagin dual $\hat G  = {\rm Hom}(G, k^{\lcross})$ of $G$. Indeed,
by Example 2 of \pa 1.2, the group ${\rm Gr}(k^G)$
is isomorphic to $\HC^1(k^G, k)$. In this situation, the identification $\HC^1(k^G, k) \cong H^1(G, k^{\lcross})$ given by Theorem 1.5 is in fact an isomorphism of groups,
and one sees that   $H^1(G, k^{\lcross})$
is isomorphic to $\hat G$.

\medskip

\noindent 2) For any finite subgroup $G$ of a group $L$, the group ring $k[L]$ is canonically equipped with a {\hbox{$k^G$-comodule}} algebra structure $\Delta _{k[L]}$ given by
$\displaystyle \Delta _{k[L]}(h) = \sum _{g \in G}ghg^{-1} \ot \delta _g$, for any $h \in L$, and extended by linearity. Theorem 1.5 claims the isomorphism $\HC^*(k^G, k[L]) \cong \Hr^*(G, k[L]^{\lcross})$.
However  the computation of the group of units in $k[L]$ is in general a very difficult problem: the group  $k[L]^{\lcross}$ is  known only for some particular groups $L$.

\goodbreak

\bigskip
\noindent {\sl 1.4. An explicit example where the Hopf algebra is not an algebra of functions on a group.} Let here~$k$ be a field and $H_4$ be the Sweedler four-dimensional Hopf algebra over $k$. Recall that 
$H_4$ is generated by two elements $g$ and $h$ submitted to the relations:
$$g^2 = 1, \quad h^2 = 0, \quad gh + hg = 0.$$ 
On the generators, the comultiplication, the antipode, and the counit of $H_4$ are given by 
$$\matrix{\Delta (g) = g \ot g, \hfill &  \quad  \Delta (h) = h \ot g + 1\ot h, \hfill \cr
\sigma (g) = g, \hfill &  \quad \sigma (gh) = gh, \hfill \cr
 \varepsilon(g) = 1, \hfill &  \quad \varepsilon(h) = 0. \hfill\cr}$$
Denote by $E_2$ the algebra of dual numbers, viewed as the subalgebra of $H_4$ generated by $h$. Via $\Delta$, the algebra $E_2$ is naturally endowed with a  structure of $H_4$-comodule algebra.

\medskip

\noindent {\tencmb Proposition 1.6.} {\sl There is an   equality of groups
$$\HC^0(H_4, E_2) = k^{\lcross}$$ and an isomorphism of pointed sets
$$ \HC^1(H_4, E_2) \cong  \{ 1\ot 1 , 1\ot g\}.$$

}

\medskip

\noindent {\sl Proof.} The proof consists in calculating explicitely the invariants on the $0$-level (we leave this point to the reader) and in writing down the cocycle relations on generic elements
on the $1$-level. The computation of $\ZC^1(H_4, E_2)$ is lightened by remarking that $E_2 \ot H_4  =  \{ 1\ot U + h \ot V \mid U , V \in  H_4\}$ and that
$H_4 = F \oplus Fh$, 
where $F$ is the sub-Hopf algebra of $H_4$ generated by $g$. 
The cocycle relation is then equivalent to the following system of two  conditions on $\Delta (U)$ and $\Delta (V)$:
$$\left\{\matrix{\Delta (U) \  = & U \ot U + Uh \ot V \hfill & (1)\cr
\Delta (V) \ = & (Ug) \ot V + V \ot U + (Vh) \ot V. \hfill & (2) \cr}\right.$$ 
In Equation (1), if one replaces $U$ by $x + yh$ and $V$ by $z + th$, with $x,y, z, t \in F$, one gets a system of four equations in $x,y, z, t$.
Solving them, one deduces  $U$ and $V$, which automatically satisfy Equation~(2). 
\goodbreak 

\noindent Finally one obtains
$$ \ZC^1(H_4, E_2) = \{ X_u , Y_u \mid u  \in k\},$$ 
where the elements $X_u$ and $Y_u$ of $\ZC^1(H_4, E_2)$ are given by
$$\matrix{X_u \hfill& = \ 1 \ot 1 + u(1\ot h) - u(h \ot 1) + u(h\ot g)  - u^2(h\ot h)\hfill \cr
Y_u \hfill& = \ 1 \ot g + u(1\ot gh) - u(h \ot g) + u(h\ot 1)  - u^2(h\ot gh). \hfill \cr}$$
The distinguished point of $ \ZC^1(H_4, E_2) $ is $X_0 = 1 \ot 1$. 
One may observe that $\ZC^1(H_4, E_2)$  contains
a group, the set $\{ X_u \mid u  \in k\}$, which acts on the right on $\ZC^1(H_4, E_2)$ by way of the multiplication in $(E_2 \ot H_4)^{\lcross}$. Indeed, for any $u, v \in k$ one has
the formulae:
$$X_u X_v = X_{u+v}, \quad Y_u X_v = Y_{u+v}.$$ 

It remains to describe the action of $E_2^{\lcross}$ on $\ZC^1(H_4, E_2)$. A generic element in $E_2^{\lcross}$ is of the form $\alpha + \beta h$, with $\alpha \in k^{\lcross}$ and $\beta \in k$.
A direct computation gives the identities
$$ X_u \leftharpoonup (\alpha + \beta h) = X_{u+\beta/\alpha} \quad {\hbox{\rm  and}} \quad Y_u \leftharpoonup (\alpha + \beta h) = Y_{u+\beta/\alpha},  $$
from which we deduce the isomorphism  $ \HC^1(H_4, E_2) \cong  \{ X_0 , Y_0\} = \{ 1\ot 1 , 1\ot g\}.$ \dm 

\goodbreak

\bigskip
\noindent {\sl 1.5. Deforming the Hopf module structure with a cocycle.} Let $H$ be a Hopf algebra and $E$ be an  {\hbox{$H$-comodule}} algebra. We show how 
the natural  structure of  $(H,E)$-Hopf module on $E$ may be deformed with the help
of a Hopf cocycle. To this end, for any element $X$ of $E \ot H$, denote by $\Delta _E^X$ the map from $E$ to $E \ot H$ given on $x \in E$ by
$$\Delta _E^X(x) = X\Delta _E(x).$$  One has then the following result:

\medskip

\goodbreak

\noindent {\tencmb Proposition 1.7.} {\sl Let $H$ be a Hopf algebra, $E$ be an  $H$-comodule algebra, and  $X$ be an element of $(E\ot H)^{\lcross}$. Then

\itemitem{1)}  the element $X$ is a Hopf $1$-cocycle if and only if  $(E, \Delta _E^X)$ is an  $(H,E)$-Hopf module; 

\itemitem{2)} two Hopf $1$-cocycles $X$ and $X'$ are  cohomologous if and only if the $(H,E)$-Hopf modules $(E, \Delta _E^X)$ and $(E, \Delta _E^{X'})$ are isomorphic.

}

\goodbreak
\medskip

\noindent {\sl Proof.} 1) Let us prove that $\Delta _E^X$ defines a coaction on $E$ if and only if $X$ belongs to $\ZC ^1(H, E)$. 
Suppose that  $X$ is a Hopf $1$-cocycle.  We have to show the two identities $(\Delta _E^X \ot \id _H)\circ \Delta _E^X = (\id _E \ot \Delta _H)\circ \Delta _E^X$
and $(\id _E \ot \varepsilon_H)\circ \Delta _E^X = \id _E$. Pick an element $x$ in $E$. On the one hand, since $\Delta _E$ is a morphism of algebras, one has the equalities
$$\eqalign{\bigl((\Delta _E^X \ot \id _H)\circ \Delta _E^X\bigr) (x) & = (\Delta _E^X \ot \id _H) \bigl(X \Delta _E (x)\bigr)\cr
 & = (X \ot 1)\Bigl((\Delta _E \ot \id _H) \bigl(X \Delta _E (x)\bigr)\Bigr) \cr
& = \Bigl((X\ot 1)\bigl( (\Delta _E \ot \id _H) (X)\bigr) \Bigr)\Bigl(\bigl((\Delta _E \ot \id _H)\circ \Delta _E\bigr)(x)\Bigr). \cr}$$
On the other hand, the following equalities hold:
$$\eqalign{\bigl((\id _E \ot \Delta _H)\circ \Delta _E^X\bigr) (x) & = (\id _E \ot \Delta _H) \bigl(X \Delta _E (x)\bigr)\cr
& = \bigl((\id _E \ot \Delta _H)(X)\bigr)\Bigl(\bigl((\id _E \ot \Delta _H)\circ \Delta _E\bigr) (x)\Bigr). \cr}$$
Since $\bigl((\Delta _E \ot \id _H)\circ \Delta _E\bigr)(x)$ is equal to $ \bigl((\id _E \ot \Delta _H)\circ \Delta _E\bigr) (x)$, it remains to remark that the identity 
$(X\ot 1)\bigl( (\Delta _E \ot \id _H) (X)\bigr) = (\id _E \ot \Delta _H)(X)$ is exactly the cocycle relation $d^2(X)d^0(X) = d^1(X)$.
In a similar way,  using Remark 1.3(a) and the identity
$(\id _E \ot \varepsilon_H)\circ \Delta _E = \id _E$, one proves the equality $(\id _E \ot \varepsilon_H)\circ \Delta _E^X = \id _E$.

The map $\Delta _E$ is a morphism of algebras, whence for any $x$ and $x'$ in $E$, one has
the equality $\Delta _E(xx') = \Delta _E(x)\Delta _E(x')$. So one gets $X\Delta _E(xx') = X\Delta _E(x)\Delta _E(x')$, or 
$\Delta _E^X(xx') = \Delta _E^X(x)\Delta _E(x')$.  This proves that $(E, \Delta _E^X)$ is an  $(H,E)$-Hopf module, where
 the $E$-module structure of $E$ is still given by the multiplication.

Conversely, assume that  $\Delta _E^X$ endows $E$ with a structure of $(H,E)$-Hopf module. Applying the 
identity $(\Delta _E^X \ot \id _H)\circ \Delta _E^X = (\id _E \ot \Delta _H)\circ \Delta _E^X$ to the element $x=1$, one 
obtains the cocycle relation for~$X$.

2) Suppose now given two cohomologous Hopf $1$-cocycles $X$ and $X'$. 
Let $x$ be an element of $E^{\lcross}$ such that $X' = (d^1x^{-1}) X (d^0x)$. One easily checks that  $\tau _x : E \lr E$, the left multiplication by $x$,
realizes an isomorphism of $(H,E)$-Hopf module from $(E, \Delta _E^{X'})$ to $(E, \Delta _E^X)$.

Conversely, assume that for two Hopf $1$-cocycles $X$ and $X'$, there exists an isomorphism of \hbox{$(H,E)$-Hopf} modules $\varphi : (E, \Delta _E^X) \lr (E, \Delta _E^{X'}) $.
By $E$-linearity, $\varphi$ is entirely determined by $\varphi (1)$, more precisely $\varphi = \tau _{\varphi (1)}$. Since $\varphi$ is surjective, the element
$\varphi (1)$ is invertible in $E$. The comodule compatibility 
relation $\Delta _E^{X'} \circ \varphi = (\varphi \ot \id _H) \circ \Delta _E^X$ then implies  $d^1 (\varphi (1))X = X'd^0 (\varphi (1))$. \dm

\bigskip
\noindent {\sl 1.6. The cohomology exact sequence associated to a sub-comodule algebra.}
  Let $H$ be a Hopf algebra. By the very definition, any  morphism $\varphi : D \lr E$ of $H$-comodule algebras gives rise to a group homomorphism 
$\HC^0(\varphi) : \HC^0(H, D) \lr \HC^0(H, E)$ and to a morphism of  pointed sets
{\hbox{$\HC^1(\varphi) : \HC^1(H, D) \lr \HC^1(H, E)$}}. Our aim is to produce an exact sequence in cohomology associated to any 
inclusion $\varphi : D \longrighthook E$ of $H$-comodule algebras. To this purpose, we state the following lemma, which is a slight generalization to the cosimplicial
case of Serre's exact sequence enounced in the framework of non-abelian cohomology theory of groups.

\medskip

\noindent {\tencmb Lemma 1.8.} {\sl Let $\varphi : \AC ^* \lr \BC ^*$ be an injective morphism of two pre-cosimplicial groups
$$ \AC ^*= \xymatrix{ A^0  \ar@<1.3ex>[r]^{d^0} 
\ar@<-1.3ex>[r]^{d^1}  & \ A^1
\ar@<2.3ex>[r]^{d^0} 
\ar@<0ex>[r]^{d^1} 
\ar@<-2.3ex>[r]^{d^2}& \ A^2} \quad \quad {\hbox{\sl and}} \quad \quad
 \BC ^*= \xymatrix{ B^0  \ar@<1.3ex>[r]^{d^0} 
\ar@<-1.3ex>[r]^{d^1}  & \ B^1
\ar@<2.3ex>[r]^{d^0} 
\ar@<0ex>[r]^{d^1} 
\ar@<-2.3ex>[r]^{d^2}&  \ B^2.}
$$
Let $\CC ^* = \BC ^*/ \varphi(\AC ^*)$  be the pre-cosimplicial left quotient object in the category of pointed sets and let $\pi : \BC ^* \lr \CC ^*$ be the quotient map. Then there is an exact sequence
of pointed sets 

$$1 \lr \HBB^0(\AC ^*)  \ \hfll{ \HBBp^0(\varphi)}{}   \     \HBB^0(\BC ^*)  \ \hfll{ \HBBp^0(\pi)}{}   \    \HBB^0(\CC ^*) 
 \ \hfll{\partial}{}   \   \HBB^1(\AC ^*)  \ \hfll{ \HBBp^1(\varphi)}{}   \    \HBB^1(\BC ^*) .$$
\goodbreak
Moreover, if $\varphi(A^i)$ is for $i = 0,1,2$ a normal subgroup of $B^i$, then the above exact sequence extends to the right in the following way:

$$1 \lr \HBB^0(\AC ^*)  \ \hfll{ \HBBp^0(\varphi)}{}   \     \HBB^0(\BC ^*)  \ \hfll{ \HBBp^0(\pi)}{}   \    \HBB^0(\CC ^*) 
 \ \hfll{\partial}{}   \   \HBB^1(\AC ^*)  \ \hfll{ \HBBp^1(\varphi)}{}   \    \HBB^1(\BC ^*)   \hfll{ \HBBp^1(\pi)}{}   \HBB^1(\CC ^*).$$

}

\Dem The connecting morphism $\partial$ is obtained by usual diagram-chasing. We leave the reader check the functoriality of $\HBB^*$ as well as the exactness of the two sequences. \dm 

\medskip

We mention here that the definition of the non-abelian $0$-cohomology object $\HBB^0({\cal A}^*)$ as an equalizer does in fact not require any algebraic structure on the set $A^0$.
This observation leads to the following definition. For any inclusion of $H$-comodule algebras {\hbox {$\varphi : D \longrighthook E$}}, we introduce the {\sl relative non-abelian $0$-cohomology set}
$$\HC^0(H, D \longrighthook E) = \HBB^0\bigl({\cal C}^{\lcross}_{\leqp2}(H, E) / {\cal C}^{\lcross}_{\leqp2}(H, D)\bigr),$$ where $ {\cal C}^{\lcross}_{\leqp2}(H, E) / {\cal C}^{\lcross}_{\leqp2}(H, D)$
is the pre-cosimplicial diagram of pointed sets
$$\xymatrix{ E^{\lcross}/D^{\lcross}  \ar@<1.3ex>[r]^{d^0\ \ \ \ \ \ \ \ } 
\ar@<-1.3ex>[r]^{d^1\ \ \ \ \ \ \ \ }  & \ (E\ot H)^{\lcross}/(D\ot H)^{\lcross}
\ar@<2.3ex>[r]^{d^0\ \ \ \ \ \ } 
\ar@<0ex>[r]^{d^1\ \ \ \ \ \ } 
\ar@<-2.3ex>[r]^{d^2\ \ \ \ \ \ }& \ (E\ot H\ot H)^{\lcross}/(D\ot H\ot H)^{\lcross}}.$$ 

In the particular case where  $D^{\lcross} $, $(D\ot H)^{\lcross}$,  and  $(D\ot H\ot H)^{\lcross}$, are  normal subgroups respectively of
$E^{\lcross} $, $(E\ot H)^{\lcross}$,  and  $(E\ot H\ot H)^{\lcross}$, then  ${\cal C}^{\lcross}_{\leqp2}(H, E) / {\cal C}^{\lcross}_{\leqp2}(H, D)$
is a pre-cosimplicial  group, and the definition $$\HC^1(H, D \longrighthook E) = \HBB^1\bigl({\cal C}^{\lcross}_{\leqp2}(H, E) / {\cal C}^{\lcross}_{\leqp2}(H, D)\bigr)$$
makes sense. Next result is a corollary of Lemma 1.8.

 \medskip

\goodbreak

\noindent {\tencmb Proposition 1.9.} {\sl Let $H$ be a Hopf algebra and  $\varphi : D \lr E$  be an injective morphism of $H$-comodule algebras. The sequence of pointed sets
$$1 \lr \HC^0(H, D)  \ \hfll{\HCp^0(\varphi)}{}   \   \HC^0(H, E)  \ \hfll{\HCp^0(\pi)}{}  \    \HC^0(H, D \longrighthook E) \  \hfll{\partial}{}   \    \HC^1(H, D)  \  \hfll{\HCp^1(\varphi)}{}  \   \HC^1(H, E) $$
is exact.
Moreover, if $D^{\lcross} $, $(D\ot H)^{\lcross}$,  and  $(D\ot H\ot H)^{\lcross}$ are normal subgroups respectively of
$E^{\lcross} $, $(E\ot H)^{\lcross}$,  and  $(E\ot H\ot H)^{\lcross}$, then the above exact sequence can be extended to the right in the following way:
$$\eqalign{1 \lr &\HC^0(H, D) \ \hfll{\HCp^0(\varphi)}{}  \   \HC^0(H, E) \ \hfll{\HCp^0(\pi)}{}   \  \HC^0(H, D \longrighthook E)  \ \hfll{\partial}{}  \   \HCp^1(H, D)  \ \hfll{\HCp^1(\varphi)}{} \cr
& \cr
& \hfll{\HCp^1(\varphi)}{} \ 
  \HC^1(H, E)  \ 
\hfll{\HCp^1(\pi)}{}  \  \HC^1(H, D \longrighthook E).\cr}$$

}

\bigskip

\goodbreak
\noindent {\tencmb 2. Links between general and restricted non-abelian Hopf cohomology theory.}
\smallskip 
\noindent In this section, $H$ is a Hopf algebra, $S$ is an $H$-comodule algebra, and $M$ is an $(H,S)$-Hopf module. 
In [8],
we introduced a cohomology theory $\Hr^*(H, M)$ that we qualify from now on as {\sl restricted}.
We first  breafly recall its definition and then compare it to our general cohomology theory under some lax technical conditions.

\goodbreak
\medskip
\noindent {\sl 2.1. Reminder on restricted non-abelian Hopf cohomology theory.}
As in [8], we endow the set $\Wr_k^n(M) = {\Hom}_k(M , M \ot H^{\otimes n})$ with a $k$-algebra structure thanks to 
the composition-type product $${\hbox {$\lodot: \Wr_k^n(M) \cross \Wr_k^n(M) \lr \Wr_k^n(M)$}}$$
given by
$$\left\{\eqalign{\varphi \lodot \varphi' & = \varphi \circ \varphi' \ \ \ {\rm  if} \ \ n = 0\cr
\varphi \lodot \varphi' & = (\id _M \ot \mu _H^{\otimes n})\circ (\id _M \ot \chi_n) \circ (\varphi \ot \id _H^{\otimes n}) 
\circ \varphi' \ \ {\rm  if} \ \ n > 0\cr}\right.$$ 
for $\varphi, \varphi' \in \Wr_k^n(M)$;
here $\chi_n: H^{\otimes n}\ot H^{\otimes n} \lr (H\ot H)^{\otimes n}$ denotes the intertwining operator defined by
$$\chi_n\bigl((a_1 \ot \ldots \ot a_n) \ot (b_1 \ot \ldots \ot b_n)\bigr) = (a_1 \ot b_1) \ot \ldots \ot (a_n \ot b_n).$$
Denote by  $\Wr_S^n(M)$ the subalgebra ${\Hom}_S(M , M \ot H^{\otimes n})$ of $\Wr_k^n(M)$, where the 
$S$-module structure on $M \ot H^{\otimes n}$ is given by
$(m \ot \underline {h}) s = ms \ot \underline {h},$ for any $m \in M$, $\underline {h} \in H^{\otimes n}$, and
$s \in S$. 

Let $R$ be either the ground ring $k$ or the algebra $S$. The 
algebras $ \Wr_R^0(M)$, $ \Wr_R^1(M)$ and $ \Wr_R^2(M)$ may be organized in a pre-cosimplicial 
diagram of monoids [8, Lemma 1.1]:
$$ {\cal W}_{\leqp2}(H, M)_R = \Bigl(\xymatrix{ \Wr_R^0(M)  \ar@<1.3ex>[r]^{b^0} 
\ar@<-1.3ex>[r]^{b^1}  &  \Wr_R^1(M) 
\ar@<2.3ex>[r]^{b^0} 
\ar@<0ex>[r]^{b^1} 
\ar@<-2.3ex>[r]^{b^2}&  \Wr_R^2(M)\Bigr).}$$
The two maps $b^i: \Wr_R^0(M) \lr \Wr_R^1(M)$
($i= 0,1$)
and the three maps $b^i: \Wr_R^1(M) \lr \Wr_R^2(M)$ ($i=0,1,2$) 
are given for $\varphi\in \Wr_R^0(M)$ and $\Phi \in \Wr_R^1(M)$ by the formulae 
$$\eqalign{ b^0\varphi   & = (\id_{M} \ot \mu _H) \circ (\Delta _{M} \ot \id _H)
\circ (\varphi \ot \sigma _H) \circ \Delta _M \cr 
b^1\varphi & = (\id _{M} \ot \eta _H)\circ \varphi\cr
 b^0\Phi & = (\id _{M} \ot \mu _H \ot \id _H)
 \circ (\Delta _{M} \ot T) \circ (\Phi  \ot \sigma _H) \circ \Delta _{M}\cr
b^1\Phi & = (\id _{M} \ot  \Delta _H) \circ \Phi \cr
b^2\Phi & = (\id _{M} \ot \id _H \ot \eta _H) \circ \Phi  = \Phi \ot \eta _H, \cr
}$$
where    $T$ denotes the flip of $H \ot H$
({\sl i.e.} the automorphism of $H \ot H$ which sends an indecomposable tensor $h \ot h'$ to
$h'\ot h$).

\goodbreak
\bigskip

Recall the definitions stated in [8]. The  restricted $0$-cohomology group $\Hr^0(H, M)$ is the equalizer \hbox{$\{ \varphi \in \Aut _S(M) \ \vert \ b^1\varphi = b^0\varphi \}$}
of the pair $(b^0, b^1)$. The restricted $1$-cohomology set $\Hr^1(H, M)$ is the quotient set
$\Aut _S(M) \backslash \Zr^1(H, M)$ of the set $\Zr^1(H, M)$ of restricted Hopf $1$-cocycles of $H$ with coefficients in $M$
under the right action of the group $\Aut _S(M)$. Recall that $\Zr^1(H, M)$ is the subgroup 
$$\Zr^1(H, M) = \left\{ \Phi \in  \Wr_k^1(M)\quad \left\vert 
\eqalign{
\quad & (\Zr\Cr_1) \ \ \ \Phi(ms) = \Phi(m)s {\hbox {, for all}} \ m \in M \ {\hbox {and}} \ s \in S\hfill \cr
 &  (\Zr\Cr_2) \ \ \ (\id _M \ot \varepsilon _H) \circ \Phi = \id _M\hfill\cr
& (\Zr\Cr_3)  \ \ \ \ b^2\Phi \lodot b^0\Phi = b^1\Phi\hfill\cr}\right. \right\}$$ of $\Wr_S^1(M)$ and an element $f \in  \Aut _S(M)$
acts on the right on an element $\Phi \in \Zr^1(H, M)$  by $$(\Phi \leftharpoonup f) = b^1f^{-1}\lodot\Phi\lodot b^0f.$$

\goodbreak

\noindent We now give a new alternative description of $\Zr^1(H, M)$ which we shall need in the sequel. 

\medskip  
\noindent {\tencmb Proposition 2.1:} {\sl 
The set $\Zr^1(H, M) $ may be written as $$\Zr^1(H, M) = \{ \Phi \in  \Wr_S^1(M)^{\lcross} \quad \left\vert 
\quad  b^2\Phi \lodot b^0\Phi = b^1\Phi \}. \right. $$

}

\medskip

\Dem Let $\Phi $ be an element of $\Zr^1(H, M)$. First observe that Condition $(\Zr\Cr_1)$ means exactly that $\Phi $ belongs to $\Wr_S^1(M)$.
It suffices to prove that Condition $(\Zr\Cr_2)$ is equivalent to the {$\lodot$}-invertibility of $\Phi$ under  Condition $(\Zr\Cr_3)$. 
Set $F = \Phi \lodot \Delta _M$.   In the proof of Theorem 3.1 in [8], we showed that $\Phi$ satisfies  $(\Zr\Cr_2)$ if and only if $F$ satisfies Condition $(\Cr\Cr_2)$, that is $(\id _M \ot \eps _H) \circ F = \id _M$. 
Similarly $\Phi$ fulfils  $(\Zr\Cr_3)$ if and only if $F$ fulfils Condition $(\Cr\Cr_3)$, that is $ (F \ot \id _H) \circ F = (\id _M \ot \Delta _H) \circ F$.

1) Suppose that $\Phi$ is invertible in $\Wr_S^1(M)$ with inverse $\Phi'$. Since the comultiplication map $\Delta _M$ is invertible in $\Wr_k^1(M)$ with inverse $\Delta '_M = (\id _M \ot \sigma_H)\circ \Delta _M$,
the map $F$ is invertible in $\Wr_k^1(M)$ with inverse $F' = \Delta '_M \lodot \Phi'$.
Compose both terms of the equality $(\Cr\Cr_3)$ on the left with the map $\id _M \ot \id _H \ot \varepsilon _H$. One gets
$F = F \circ \bigl( (\id _M \ot \varepsilon _H) \circ F\bigr)$, which is equivalent to the relation
$F = F \lodot \Bigl(\bigl( (\id _M \ot \varepsilon _H){\circ F}\bigr) \ot \eta _H\Bigr)$. One may simplify by $F$, and one gets
$\bigl((\id _M \ot \varepsilon _H){\circ F}\bigr) \ot \eta _H = \id _{{\hbox{\sevenrm W}}_k^1(M)} = \id _M \otimes \eta _H$.
Applying now $\id _M \ot \varepsilon _H$ on the right, one obtains $(\Cr\Cr_2)$.

2) Conversely, assume that Condition $(\Cr\Cr_2)$ holds. We shall show that the map $F'$ defined by $F' = (\id _M \ot \sigma _H) \circ F$
is the inverse of $F$ in $\Wr_k^1(M)$. We apply therefore  $\id _M \ot \bigl(\mu _H \circ (\id _H \ot \sigma _H)\bigr)$, respectively 
 $\id _M \ot \bigl(\mu _H \circ (\sigma _H \ot \id _H)\bigr)$,
on the left to the equality $(\Cr\Cr_3)$. We get $\bigl( (\id _M \ot \varepsilon _H){\circ F}\bigr) \ot \eta _H = F \lodot F'$, respectively $\bigl( (\id _M \ot \varepsilon _H){\circ F}\bigr)  = F' \lodot F$.
By Condition $(\Cr\Cr_2)$, this exactly means that $F'$ is the inverse of $F$. 
The map $\Phi$ is therefore invertible in $\Wr_k^1(M)$ with inverse $$\Phi ' = \Delta _M \lodot F' = \Delta _M \lodot \bigl((\id _M \ot \sigma _H) \circ ( \Phi \lodot \Delta _M)\bigr).$$
It remains to show that $\Phi '$ is $S$-linear.
For  $m \in M$, we  denote the tensor
$\Phi (m) \in M \ot H$ by $m_{[0]} \ot m_{[1]}$. We have
$$\Phi' (m) = \bigl((m_0)_{[0]}\bigr)_0\ot \bigl((m_0)_{[0]}\bigr)_1\sigma_H\bigl( (m_0)_{[1]}m_1\bigr).$$

\goodbreak
For any $s \in S$, we obtain
$$\eqalign{\Phi' (ms) & = \bigl((m_0s_0)_{[0]}\bigr)_0\ot \bigl((m_0s_0)_{[0]}\bigr)_1\sigma_H\bigl( (m_0s_0)_{[1]}m_1s_1\bigr)\cr
& = \bigl((m_0)_{[0]}s_0\bigr)_0\ot \bigl((m_0)_{[0]}s_0\bigr)_1\sigma_H(s_1) \sigma_H\bigl((m_0)_{[1]}m_1\bigr)\cr
& = \bigl((m_0)_{[0]}\bigr)_0s_0\ot \bigl((m_0)_{[0]}\bigr)_1s_1\sigma_H(s_2) \sigma_H\bigl((m_0)_{[1]}m_1\bigr)\cr
& = \bigl((m_0)_{[0]}\bigr)_0s_0\varepsilon_H(s_1)\ot \bigl((m_0)_{[0]}\bigr)_1\sigma_H\bigl((m_0)_{[1]}m_1\bigr)\cr
& = \Phi '(m)s.}$$
This computation, which proves the $S$-linearity of $\Phi '$,  uses the Hopf algebra yoga. Moreover
the first and the third equalities come from  $\Delta _M(ms) = \Delta _M(m)\Delta _S(s)$, whereas the second one is a consequence of the $S$-linearity of $\Phi$. ${}$ \dm

\medskip 
Denote by $$ {\cal W}^{\lcross}_{\leqp2}(H, M) = \Bigl(\xymatrix{ \Wr_S^0(M)^{\lcross}  \ar@<1.3ex>[r]^{b^0} 
\ar@<-1.3ex>[r]^{b^1}  &  \Wr_S^1(M) ^{\lcross}
\ar@<2.3ex>[r]^{b^0 \ } 
\ar@<0ex>[r]^{b^1 \ } 
\ar@<-2.3ex>[r]^{b^2 \ }&  \Wr_S^2(M)^{\lcross}\Bigr),}$$ 
the pre-cosimplicial diagram of groups obtained by taking the $\lodot$-invertible elements of $\Wr_S^*(M)$.
Proposition 2.1 leads us to state the following result:
\medskip
\noindent {\tencmb Theorem 2.2:} {\sl Let $H$ be a Hopf algebra, $S$ be an $H$-comodule algebra, and $M$ be an $(H,S)$-Hopf module.  
One has the equality $$\Hr^*(H, M) = \HBB^*({\cal W}^{\lcross}_{\leqp2}(H, M)). $$
}

\medskip
\noindent {\sl 2.2.  Technical conditions.} In this paragraph, we first point out technical conditions
we shall need in the sequel in order to compare the general and the restricted non-abelian Hopf cohomology theories.
We show then that these conditions are fulfilled in two natural cases.
\medskip

For any $n \geq 0$, consider the linear map  $$\omega _n : \End _S(M)\ot H^{\otimes n} \lr \Wr_S^n(M) = {\Hom}_S(M , M \ot H^{\otimes n})$$
given on an undecomposable tensor $f \ot \underline {h} \in \End _S(M) \ot H^{\otimes n}$  by
$$\omega _n(f \ot \underline {h})(m) = f(m) \ot \underline {h},$$
where $m \in M$. Notice that $\omega _0 $ is the identity map of $\End _S(M)$ and
that, for any $n \geq 0$, the map $\omega _n $ is a morphism of  algebras.
For $n\geq 0$, we consider the following condition. 
\bigskip
\centerline{{\sl Condition $({\cal F}_n)$}: the map $\omega _n$ is an isomorphism of algebras.}

 \bigskip
By the very definitions, Condition $({\cal F}_0)$ always holds.
The first natural case where Condition $({\cal F}_n)$ is satisfied for all $n \geq 0$ appears when $H$ is a finitely generated free $k$-module.
We develop now a second case.

 \medskip
 Let $M^* = \Hom _k(M, k)$ be the linear dual of the $k$-module $M$. Consider  the evaluation map {\hbox {$d_M : M^* \ot M \lr k$}} given by
$d_M (\nu \ot m) = \nu (m)$, for any $\nu \in M^*$ and  $m \in  M$.

\medskip  
\noindent {\tencmb Proposition 2.3:} {\sl Condition $({\cal F}_n)$ is satisfied for all $n \geq 0$ if the two following statements both hold:

1) the Hopf algebra $H$ is free as a $k$-module or $S$ is equal to the ground ring $k$;

2) there exists a map $b_M : k \lr M \ot M^*$, called birth-map, such that

$$(\id _M \ot d_M)\circ (b_M \ot \id _M) = \id _M \ \quad {\hbox {\sl and}} \ \quad (d_M \ot \id _{M^*} )\circ ( \id _{M^*} \ot b_M) = \id _{M^*}.$$ }

By convention, we set $\displaystyle b_M (1) = \sum _ie_i \ot e^i$. With this notation, the previous two equalities are equivalent to 
$$\sum _ie_ie^i(m) = m \ \quad {\hbox {\rm  and}} \ \quad \sum _i\nu (e_i)e^i = \nu,$$
for any $m \in M$ and  $\nu \in M^*$.

\medskip

\noindent {\sl Example:}  When $M$ is a finitely generated free $k$-module
with basis $(e_j)_{j = 1, \ldots , n}$  such a birth-map $b_M$ exists and is given by
$\displaystyle b_M (1) = \sum _{j = 1}^ne_j \ot e_j^*$. Here  $(e_j^*)_{j = 1, \ldots , n}$ is the dual basis 
of  $(e_j)_{j = 1, \ldots , n}$.

\medskip

The data of a module together with an evaluation map and a birth-map abstracts the notion of duality in tensor categories 
(see [2]).

 \medskip
 \noindent {\sl Proof of Proposition 2.3:} 
First of all, we endow $M^*$ with the left $S$-module structure given by $$(s\nu )(m) = \nu (ms)$$
with $\nu \in M^*$, $m\in M$, and $s \in S$. The module $M \ot M^* \ot  H^{\otimes n}$ becomes an algebra 
through the multiplication given on two elements $m\ot \nu \ot \underline {h}$ and $m'\ot \nu' \ot \underline {h'}$ of $M \ot M^* \ot  H^{\otimes n}$ by the formula
$(m\ot \nu \ot \underline {h})(m'\ot \nu' \ot \underline {h'}) = \nu(m')m\ot \nu \ot \underline {h} \ \underline {h'}$.
We introduce the subalgebra $\Er_S^n(M)$ of $M \ot M^* \ot  H^{\otimes n}$ consisting of the
elements $m\ot \nu \ot \underline {h}$ such that, for any $s\in S$, one has  $ms\ot \nu \ot \underline {h} = m\ot s\nu \ot \underline {h}$.
Notice that under the first statement, one has the equality
$$\Er_S^n(M) = \Er_S^0(M)\ot  H^{\otimes n}.$$

\goodbreak
We show now that, under the second statement, $\Er_S^n(M)$ is isomorphic to $\Wr_S^n(M)$ as an algebra. First observe that 
the existence of a birth-map allows to  write the action of $s $ on $\nu$ as
{\hbox{$s \nu = \sum_i \nu (e_i s) e^i$}}. Moreover one has  $\sum _i e_is \ot e^i = \sum _i e_i \ot se^i,$ in other words, $b_M(1)$ belongs to~$\Er_S^0(M)$.

Consider the morphism $\lambda _n : \Er_S^n(M) \lr \Wr_S^n(M)$
defined by
$$\bigl(\lambda _n(m \ot \nu \ot \underline {h})\bigr)(m') = \nu(m')m\ot \underline {h},$$
with $m, m' \in M$,  $\nu \in M^*$, and $\underline {h} \in H^{\otimes n}$. One  checks that $\lambda _n$ is well-defined with respect to the {\hbox {$S$-invariance}} and that 
it is a morphism of algebras. 
We prove now that
under the existence of a birth-map, $\lambda _n$ is a bijection. Let us explicit the inverse map. Denote by
$\lambda '_n: \Wr_S^n(M) \lr M \ot M^* \ot  H^{\otimes n}$ the map given on an element $\Phi \in \Wr_S^n(M)$ by 
$$\lambda ' _n (\Phi) = \sum _i \Phi(e_i)_0\ot e^i \ot \Phi(e_i)_1,$$
where, for any $m \in M$, we set $\Phi(m) = \Phi(m)_0 \ot \Phi(m)_1 \in M \ot H^{\otimes n}$.

\goodbreak

The map $\lambda '_n$ takes its values in $\Er_S^n(M)$. Indeed, using the $S$-linearity of $\Phi \in \Wr_S^n(M)$ and the fact that
$b_M(1)$ belongs to $\Er_S^0(M)$, we have, for any $s \in S$:

$$\sum _i \Phi(e_i)_0\ot se^i \ot \Phi(e_i)_1 = \sum _i \Phi(e_is)_0\ot e^i \ot \Phi(e_is)_1 =  \sum _i \Phi(e_i)_0s\ot e^i \ot \Phi(e_i)_1.$$

Moreover the map $\lambda '_n$ is a morphism of algebras: for $\Phi, \Psi \in \Wr_S^n(M)$, one has

$$\eqalign{\lambda '_n (\Phi) \lambda '_n (\Psi) & = \sum _{i, j} e^i\bigl(\Psi(e_j)_0\bigr)\Phi(e_i)_0\ot e^j \ot \Phi(e_i)_1\Psi(e_j)_1 \cr
& = \sum _{i, j} \Phi(e^i\bigl(\Psi(e_j)_0\bigr)e_i)_0\ot e^j \ot \Phi(e^i\bigl(\Psi(e_j)_0\bigr)e_i)_1\Psi(e_j)_1 \cr
& = \sum _{j} \Phi( \Psi(e_j)_0)_0\ot e^j \ot \Phi(\Psi(e_j)_0)_1\Psi(e_j)_1 \cr
& = \sum _j (\Phi\lodot \Psi)(e_j)_0\ot e^j \ot (\Phi\lodot \Psi)(e_j)_1\cr
&= \lambda '_n (\Phi\lodot \Psi)}.$$

It remains to compute the two compositions $\lambda _n\circ \lambda '_n$ and $\lambda '_n\circ \lambda _n$.
One has, for any $\Phi \in  \Wr_S^n(M)$ and $m\in M$:
$$\lambda _n(\lambda '_n(\Phi))(m) = \sum _i e^i(m)\Phi(e_i)_0\ot \Phi(e_i)_1 = \Phi(\sum _i e^i(m)e_i) = \Phi(m).$$ 
On the other hand, for $m\ot \nu \ot \underline {h} \in M \ot M^* \ot  H^{\otimes n}$, one obtains
$$\lambda '_n(\lambda _n(m\ot \nu \ot \underline {h})) = \sum _i \nu(e_i)m \ot e^i  \ot \underline {h} = m \ot  (\sum _i \nu(e_i)e^i)  \ot \underline {h} = m \ot  \nu \ot \underline {h}.$$
To end the proof, we write down the following sequence of isomorphisms, the composition of which is~$\omega _n$:
$$\End _S(M)\ot H^{\otimes n} = \Wr_S^0(M)\ot H^{\otimes n} \hfll{ \lambda _0 \otimes \idr _{H}^{\otimes n}}{}   \ 
\Er_S^0(M)\ot H^{\otimes n} = \Er_S^n(M) \hfll{ \lambda _n^{-1}}{}   \Wr_S^n(M). $$
Hence $\omega _n$ is an isomorphism of algebras, {\sl i.e.} Condition $({\cal F}_n)$ is fulfilled. \dm

\goodbreak

\bigskip
\noindent {\sl 2.3. An $H$-comodule structure on $\End _S(M)$.}
Suppose from now on that Condition $({\cal F}_n)$ is satisfied for  $0 \leq n \leq 2$. We define
the morphism $\Delta _{\End _S(M)} : \End _S(M) \lr \End _S(M) \ot H$ to be the composition map 
$$\End _S(M) = \Wr_S^0(M) \hfll{b^0}{}   \ 
\Wr_S^1(M) \hfll{ \omega _1 ^{-1}}{}  \End _S(M) \ot H. $$

\eject
\medskip  
\noindent {\tencmb Lemma 2.4:} {\sl The map $\Delta _{\End _S(M)}$ endows $\End _S(M)$ with a structure of $H$-comodule algebra.}

 \medskip
 \noindent {\sl Proof:} 
As a composition of morphisms of algebras, $\Delta _{\End _S(M)}$ is  a morphism of algebras. Let us prove that
 $\Delta _{\End _S(M)}$ is coassociative. To this end, consider the following diagram in which the upper horizontal  and the left vertical  compositions are $\Delta _{\End _S(M)}$:

$$\xymatrix{\End _S(M) \ar@{=}[r]\ar@{=}[d] & \Wr_S^0(M)\ar@{=}[ld]\ar[r]^{b^0}& \Wr_S^1(M)\ar[r]^{\omega _1 ^{-1} \quad}\ar@{=}[rd]&\End _S(M) \ot H \ar[d]^{\omega _1} &\\
\Wr_S^0(M) \ar[d]_{b^0}  &   & &    \Wr_S^1(M)\ar[d]^{b^0}  \\
\Wr_S^1(M) \ar[d]_{\omega _1 ^{-1}} \ar@{=}[rd] &   & &    \Wr_S^2(M)\ar[d]^{\omega _2 ^{-1}}\ar@{=}[ld]  \\
\End _S(M)\ot H \ar[r]_{\quad \omega _1}& \Wr_S^1(M)\ar[r]_{b^1}& \Wr_S^2(M)\ar[r]_{\omega _2 ^{-1} \quad }&\End _S(M) \ot H^{\otimes 2}  &\\
& & &}$$
The pre-cosimplicial relation $b^0b^0 = b^1b^0$ implies the commutativity of the inner octogon, hence of the whole diagram. One may see that the lower horizontal composition
is $\id _{\End _S(M)} \ot \Delta _H$ and that the right vertical composition  is $\Delta _{\End _S(M)} \ot \id _H$. This shows the coassociativity of
$\Delta _{\End _S(M)}$.

 The compatibility with the counit $(\id _{\End _S(M)}\ot \varepsilon _H)\circ \Delta _{\End _S(M)} = \id _{\End _S(M)}$ is a consequence of the relation
$(\id _{\End _S(M)}\ot \varepsilon _H)\circ b^0(\varphi) = \varphi$, which holds for all $\varphi \in \End _S(M)$. \dm

\medskip

This construction allows us to define the cohomology of the Hopf algebra $H$ with values  in the {\hbox{$H$-comodule}} algebra $\End _S(M)$. So, 
under the hypothesis that Condition $({\cal F}_n)$ is satisfied for  $0 \leq n \leq 2$, the cohomology sets $\HC^i(H, \End _S(M))$ ($i= 0,1$) make sense.

\bigskip
\noindent {\sl 2.4. The Comparison Theorem.} We are now able to compare  restricted and general non-abelian Hopf cohomology theories.

\medskip  
\noindent {\tencmb Proposition 2.5:} {\sl Let $H$ be a Hopf algebra, $S$ be an $H$-comodule algebra, and $M$ be an $(H,S)$-Hopf module such that Condition $({\cal F}_n)$ is satisfied for  $0 \leq n \leq 2$.
The pre-cosimplicial groups ${\cal C}^{\lcross}_{\leqp2}(H, M)$  and ${\cal W}^{\lcross}_{\leqp2}(H, \End _S(M))$ are isomophic. }
\medskip

\Dem The map $\omega _n$ is an isomorphism of algebras since Condition $({\cal F}_n)$ holds. Moreover, one checks the equalities
$$\omega _j d^i = b^i \omega _{j-1}$$
for any $1 \leq j \leq 2$ and $0 \leq i \leq j$. \dm

\medskip  

Theorem 2.2 and Proposition 2.5 imply the following result:

\medskip 
\noindent {\tencmb Theorem 2.6:} {\sl Let $H$ be a Hopf algebra, $S$ be an $H$-comodule algebra, and $M$ be an $(H,S)$-Hopf module such that Condition $({\cal F}_n)$ is satisfied for  $0 \leq n \leq 2$.
Then there is an  equality of groups
$$\HC^0(H, \End _S(M)) = \Hr^0(H, M)$$ and an isomorphism of pointed sets
$$ \HC^1(H, \End _S(M))\cong \Hr^1(H, M).$$ }

\bigskip

\noindent {\tencmb 3. Hopf torsors.}
\smallskip 
\noindent In this section, we define Hopf torsors. They generalize the classical torsors used in the framework of groups. 
We show that Hopf torsors are classified by a general non-abelian Hopf $1$-cohomology set.

\medskip
\noindent {\sl 3.1. Definition of Hopf torsors.}
 Let $E$ be an algebra and  $T$ be a left $E$-module. For any $u \in T$, consider the $E$-linear map $\vartheta _u : E \lr T$ defined on $x\in E$ by $\vartheta _u(x) = ux$. 
Denote by $T^{\lcross}$  the set $$T^{\lcross} = \{u \in T \mid \vartheta _u {\hbox{\rm  \ is bijective}} \}. $$

\smallskip

\noindent  From now on, we deal with  $E$-modules $T$ such that $T^{\lcross}$  is not empty.
For example, if $T$ is $E$ itself the above set coincides with the group $E^{\lcross}$ 
of invertible elements of the algebra $E$. 
Moreover, for any $E$-module $T$, observe that $T^{\lcross}$ inherits the structure of an $E^{\lcross}$-set.
In the following lemma, we collect several technical results about $T^{\lcross}$.

\medskip

\noindent  {\tencmb Lemma 3.1:} {\sl Let $E$ be an algebra and  $T$ be a left $E$-module such that the set $T^{\lcross}$ is not empty. 

 \noindent 1) Let $u$ be  an element of  $T^{\lcross}$. Then 
 $\vartheta _u^{-1}(v)$ is, for any $v \in T$, the unique element of $E$ such that $v\vartheta _u^{-1}(v) = v$. 

\noindent 2) Let $v$ and $v'$ be two elements in $T$ and  $u$ be an element in  $T^{\lcross}$. Then one has the identity
$\vartheta _u^{-1}(v)\vartheta _u^{-1}(v') = \vartheta _u^{-1}(v\vartheta _u^{-1}(v'))$.

\noindent 3) For any $u \in T^{\lcross}$, the map $\vartheta _u$ realizes a bijection between 
$E^{\lcross}$ and $T^{\lcross}$. 

}

\medskip
\Dem The first point is a direct consequence of the definition of 
 $\vartheta _u$. To show the second point, one writes  $u\vartheta _u^{-1}(v)\vartheta _u^{-1}(v') = v\vartheta _u^{-1}(v') = u\vartheta _u^{-1}(v\vartheta _u^{-1}(v'))$, and  concludes by uniqueness.
Let us prove the third point. We have to show that, for any $u \in T^{\lcross}$, the bijection  $\vartheta _u : E \lr T$ restricts to a bijection between 
$E^{\lcross}$ and $T^{\lcross}$. For any $u \in T^{\lcross}$, the set $\vartheta _u (E^{\lcross})$ is contained in $T^{\lcross}$. Indeed if  $x$ belongs to $E^{\lcross}$, one has
$\vartheta _{ux} = \vartheta _u \circ \tau _x$, where $\tau _x $ denotes the left multiplication by $x$, which is bijective. 
The induced map remains injective.  To prove that it is surjective, it is sufficient to show that   $\vartheta _u^{-1}(v)$ belongs to $E^{\lcross}$ for any $v \in T^{\lcross}$. By point 2), one has
$v\vartheta _u^{-1}(v)\vartheta _v^{-1}(u) = v\vartheta _v^{-1}(u) = u$, so $\vartheta _u^{-1}(v)\vartheta _v^{-1}(u) = 1$. 
\dm

\medskip
Let $H$ be a Hopf algebra, $E$ be an $H$-comodule algebra, and $(T, \Delta _T)$ be an  $(H, E)$-Hopf module. 
In this situation, the tensor product $T \ot H$ is an $E\ot H$-module and $(T \ot H)^{\lcross}$ makes sense. 
Notice that if $u $ belongs to $T^{\lcross}$,  then $u \ot 1$ lies in $(T\ot H)^{\lcross}$, since $\vartheta _{u \otimes 1} = \vartheta _u \ot \id _H$.
In particular, if $T^{\lcross}$ is non-empty, so is $(T\ot H)^{\lcross}$.

We introduce now  the set $$T^\bullet = \{u \in T^{\lcross} \  \vert \ \Delta _T(u) \in (T\ot H)^{\lcross}\}.$$

\bigskip

\goodbreak
\noindent  {\sl Definition 3.2:} Let $H$ be a Hopf algebra, $E$ be an $H$-comodule algebra. An $(H, E)$-{\sl Hopf torsor} is an $(H, E)$-Hopf module $(T, \Delta _T)$
such that the set $T^\bullet$ is non-empty.
\medskip
In particular $E$ is an $(H, E)$-Hopf torsor. Indeed $E$ is an $(H, E)$-Hopf module and  $\Delta _E$ being a morphism of algebras, $\Delta _E$ sends any element of $E^{\lcross}$ into $(E\ot H)^{\lcross}$. 

\medskip
We denote by $\tors (H,E)$ the set of $(H, E)$-Hopf torsors. It is pointed with distinguished point $(E, \Delta _E)$.
Two $(H, E)$-torsors $(T, \Delta _T)$ and $(T', \Delta _{T'})$ are {\sl equivalent}
if $T$ and $T'$ are isomorphic  as $(H,E)$-Hopf modules.
We denote by $\Tors (H, E)$ the set of equivalence classes of $(H, E)$-torsors; it is pointed with distinguished point
the class of $(E, \Delta _E)$.

\medskip

\noindent  {\tencmb Lemma 3.3:} {\sl Let $T$ be an $(H, E)$-Hopf torsor. 
 Then the sets $T^\bullet$ and $T^{\lcross}$ coincide.

}

\medskip

\noindent  {\sl Proof:} 
Pick $v $ in $T^{\lcross}$ and $u$ in $T^\bullet$.
One has $v = u\vartheta _u^{-1}(v)$, thus $\Delta _T(v) = \Delta _T(u)\Delta _E(\vartheta _u^{-1}(v))$. By definition, the term $\Delta _T(u)$ belongs to $(T\ot H)^{\lcross}$, and the factor
$\Delta _E(\vartheta _u^{-1}(v))$ is invertible in $E \ot H$ since $\Delta _E$ is a morphism of algebras. In the same way as $E^{\lcross}$ acts on
$T^{\lcross}$, the group  $(E\ot H)^{\lcross}$ acts on
$(T\ot H)^{\lcross}$, hence $\Delta _T(v)$ is an element of $(T\ot H)^{\lcross}$, in other words $v$ belongs to $T^\bullet$.
\dm

\medskip

\bigskip
\noindent {\sl 3.2. The non-abelian $1$-Hopf cohomology set and Hopf torsors.}
As in the world of groups, the Hopf torsors are classified by a non-abelian $1$-cohomology set. We detail this point now.
\bigskip

\noindent  {\tencmb Theorem 3.4:} {\sl Let $H$ be a Hopf algebra and $E$ be an $H$-comodule algebra. There is an isomor\-phism of pointed sets
$$ \HC^1(H, E) \cong \Tors (H, E).$$}

\medskip

\noindent  {\sl Proof:} We construct a map $\tilde {\T}: \ZC^1(H,E) \lr \tors (H,E)$ in the following way. For any Hopf {\hbox {$1$-cocycle}}  $X$, let 
$\tilde {\T}(X)$ be the $(H,E)$-Hopf module $(E, \Delta _E^X)$ defined in \pa 1.5. It is clearly  a torsor (indeed $T^\bullet$ contains for example  the unit of $E$). By Proposition 1.7, the map
$\tilde {\T}$ induces a map {\hbox{${\T}: \HC^1(H,E) \lr \Tors (H,E)$}} on the quotients.

The injectivity of ${\T}$ is a direct consequence of Proposition 1.7. Let us prove that  ${\T}$ is surjective. Take a torsor $(T, \Delta _T)$ and $u \in T^{\bullet}$. 
By definition, $\Delta _T(u)$ belongs to $(T\ot H)^{\lcross}$. 
Applying the map $\vartheta ^{-1}_{u \otimes 1} = \vartheta ^{-1}_u \ot \id _H$, we
define the element
$$X_T = (\vartheta ^{-1}_{u \otimes 1} \circ \Delta _T) (u) = \bigl((\vartheta ^{-1}_{u} \ot \id _H) \circ \Delta _T\bigr) (u),$$  which belongs to $(E \ot H)^{\lcross}$.
Writing $\Delta _T (u) = u_0 \ot u_1$, one gets $X_T = \vartheta ^{-1}_{u}(u_0) \ot u_1$. Let us compute  the product
$(u \ot 1 \ot 1)(d^2(X_T)d^0(X_T))$. First remark that we have the equalities $$(u \ot 1 \ot 1)d^2(X_T) = (\vartheta _u\vartheta ^{-1} _u \ot \id _H) (\Delta _T(u)) \ot 1 = \Delta _T(u) \ot 1.$$
On the other hand, we  write $$d^0(X_T) = (\Delta _E \circ \vartheta ^{-1} _u \ot \id _H)(\Delta _T(u)) = \Delta _E (\vartheta ^{-1} _u (u_0))\ot u_1.$$
By multiplying the two expressions, we get
$$\eqalign{(u \ot 1 \ot 1)d^2(X_T)d^0(X_T) &= \Delta _T(u)\Delta _E (\vartheta ^{-1} _u (u_0))\ot u_1\cr & =
\Delta _T(u\vartheta ^{-1} _u (u_0))\ot u_1\cr & = \Delta _T(u_0) \ot u_1\cr & = u_0 \ot \Delta _H(u_1).}$$
Finally we obtain
$$d^2(X_T)d^0(X_T) = \vartheta ^{-1} _u(u_0) \ot \Delta _H(u_1) = (\id _E \ot \Delta _H)(X_T) = d^1(X_T).$$ Hence $X_T$ is a Hopf $1$-cocycle.

We  show now that the torsors $(T, \Delta _T)$ and $\tilde {\T} (X_T) = (E, \Delta_E^{X_T})$ are equivalent. The wished isomorphism of Hopf modules between
$(E, \Delta_E^{X_T})$ and $(T, \Delta _T)$ is given by the map $\vartheta _u$. Indeed, for any element $x \in E$, one has  the equalities
$$\eqalign {\bigl((\vartheta _u \ot \id _H)\circ \Delta_E^{X_T}\bigr)(x) &= (\vartheta _u \ot \id _H)\bigl(X_T \Delta _E(x)\bigr)\cr & =
 (\vartheta _u \ot \id _H)\Bigl(\bigl((\vartheta ^{-1}_{u} \ot \id _H) \circ \Delta _T\bigr) (u) 
\Delta _E(x)\Bigr)\cr
& = \Delta _T(u) 
\Delta _E(x)\cr
& = \Delta _T(u x)\cr
& = (\Delta _T\circ \vartheta _u)(x).\cr
}$$ \dm

\medskip

\noindent {\sl Example 3.5.} Let,  as in \pa 1.4, $H_4$ be the Sweedler four-dimensional Hopf algebra over a field  $k$
and~$E_2$ be the algebra of dual numbers. The image of $\tilde {\T}$ in  $\tors (H,E)$ consists of the $(H_4, E_2)$-modules 
$T_{X_a} = (E_2, \Delta ^{X_a})$ and $T_{Y_a} = (E_2, \Delta ^{Y_a})$, where $a$ runs through $k$ and where the coactions are explicitely given by
$$\eqalign{&\Delta ^{X_a}(1) \hfill  = X_a \hfill = \ 1 \ot 1 + a(1\ot h) - a(h \ot 1) + a(h\ot g)  - a^2(h\ot h)\hfill \cr
&\Delta ^{X_a}(h) \hfill  = X_a\Delta (h) = \ 1 \ot h + h \ot g - a(h\ot h) \hfill \cr}$$
and
$$\eqalign{&\Delta ^{Y_a}(1) \hfill  = Y_a   = \ 1 \ot g + a(1\ot gh) - a(h \ot g) + a(h\ot 1)  - a^2(h\ot gh)\cr
&\Delta ^{Y_a}(h) \hfill  = Y_a \Delta (h)    = \ 1 \ot gh + h\ot 1  - a(h\ot gh). \hfill \cr}$$

\noindent Up to isomorphism, only two equivalence classes of torsors remain: those consisting in the class of $(E_2, \Delta)$ itself
and the class of $(E_2, \Delta')$, where $$ \Delta'(1)  = 1\ot g \quad \quad{\hbox{\rm  and}} \quad\quad \Delta'(h)   = h\ot 1 + 1\ot gh.$$

\medskip

\noindent  {\sl Remark 3.6:} Suppose that the algebras $E$ and $H$ are both commutative. Let $T$ and $T'$ be two $(H, E)$-Hopf torsors. Endow $T'$ with the symmetric $E$-bimodule action.
One may easily check that the tensor product $T\ot _ET'$ 
is also an $(H, E)$-Hopf torsor with coaction given by {\hbox{$\Delta _{T\otimes _ET'} (t\ot t') = t_0 \ot t'_0 \ot t_1t'_1$.}} Indeed  the set $(T\ot _ET')^{\bullet}$ 
contains all the elements $u\ot u'$, where $u$ belongs to $T^{\bullet}$ and $u'$ to $T'^{\bullet}$. 
Whence $\tors (H, E)$ is a 
monoid with product $\ot _E$. Under these hypothesis of commutativity, we already noticed that $\ZC^1(H,E)$ and $\HC^1(H,E)$
are groups (Remark 1.3(b)).  The map $\tilde {\T}: \ZC^1(H,E) \lr \tors (H,E)$ is then a morphism of monoids. Following Theorem 3.4, the product of $\tors (H, E)$
induces a group structure on the quotient~$\Tors (H, E)$.

\bigskip

\noindent {\sl 3.3. Comparison with the group case.}
Let us show how to relate Definition 3.2 to the usual notion of torsors. Given a finite group  $G$ and a $G$-group $A$, 
a  {\sl $(G, A)$-group torsor} is a non-empty left $G$-set $P$ on which $A$ acts on the
right in a compatible way with the $G$-action and such that $P$ is an affine space over $A$ (see [10]).
Denote by ${\rm Tors}(G,A)$ the set of isomorphism classes of $(G, A)$-group torsors, which is known to be isomorphic to $\Hr^1(G,A)$ (Proposition I.33 in [10]).
If $P$ is a $(G, A)$-torsor, its class in ${\rm Tors}(G,A)$ is written $[P]$.

\medskip

\noindent  {\tencmb Proposition 3.7:} {\sl Let $G$ be a finite group, let $k^G$ be the Hopf algebra of the functions on  $G$, and $E$ be an $k^G$-comodule algebra. 
For any  $(k^G, E)$-Hopf torsor $T$, the set $T^{\lcross}$ is a  $(G, E^{\lcross})$-group torsor.
 }

\bigskip

\noindent  {\sl Proof:} As previously observed, $T^{\lcross}$ is an $E^{\lcross}$-set.
By \pa 1.3, the group $E^{\lcross}$ is equipped with a $G$-group structure.
In the same way, if one writes
$\displaystyle \Delta_T(u) =  \sum _{g \in G} {^g\!u} \ot \delta _g$ for $u \in T$, one deduces an action of the  group~$G$ on the set $T$.
By Lemma 3.3, if $u$ belongs to $ T^{\lcross}$, then the element $\Delta _T(u)$ belongs to $(T\ot k^G)^{\lcross}$, which is easily seen to be isomorphic
to $(T^{\lcross})^{|G|}$. So, for any $g \in G$, the element $^g\!u$ belongs to $T^{\lcross}$, hence $T^{\lcross}$ is a $G$-group. The compatibility of the two $G$-structures on $E^{\lcross}$ and $T^{\lcross}$
is a consequence of the $(k^G,E)$-Hopf module structure of $T$. The fact that $T^{\lcross}$ is an affine space over $E^{\lcross}$ is precisely the 
bijectivity of $\vartheta _u$ proved in Lemma 3.1 for any $u \in T^{\lcross}$. \dm 

\medskip

Denote by $c : \tors (k^G, E) \lr \Torsr(G, E^{\lcross})$ the map defined for any $(k^G, E)$-torsor $T$ by 
$$c(T) = [T^{\lcross}].$$

\medskip
\noindent {\tencmb Corollary  3.8.}
{\sl Let $G$ be a finite group and $E$ be a  $k^G$-comodule algebra.
The map $c$ induces a bijection of pointed sets
$$  \Tors (k^G, E) \cong \Torsr(G, E^{\lcross}).$$
}

\Dem The isomorphism is a direct consequence of Theorem 1.5, Theorem 3.4 of this article, and Proposition I.33 in [10]. It is given by the sequence of isomorphims
$$\Tors (k^G, E) \cong \HC^1(k^G, E) \cong \Hr^1(G, E^{\lcross}) \cong \Torsr(G, E^{\lcross}).$$
Let $T$ be a $(k^G, E)$-torsor and $u$ an element of $T^{\lcross}$. The sequence of isomorphims associates to $T$ the class of the $(G, E^{\lcross})$-torsor
$E^{\lcross}_T$ defined as follows. As a set $E^{\lcross}_T$ is nothing but  $E^{\lcross}$. It is endowed with the $G$-action given for $g \in G$ and $x \in E^{\lcross}$ by $$g \rightharpoonup x = \vartheta _u^{-1}({^g\!u}){^g\!x}.$$
One verifies that $[E^{\lcross}_T] =  [T^{\lcross}] = c(T)$ in $\Torsr(G, E^{\lcross})$ via the isomorphims $\vartheta _u : E^{\lcross} \lr T^{\lcross}$. \dm

\bigskip

\noindent {\sl 3.4. Comparison with the restricted case.}
Let $H$ be a Hopf algebra, $S$ an $H$-comodule algebra, and $M$ an $(H,S)$-Hopf module. Recall that what we called {\sl $M$-torsor} in [8] is a triple
$(X, \Delta _X, \beta)$, where $\Delta _X: X \lr X \ot H$ is a map conferring $X$ 
a structure of $(H,S)$-Hopf module and {\hbox{$\beta: M \lr X$}}is an $S$-linear isomorphism.  Here we rename this datum a {\sl restricted  $M$-torsor}. 
The set of restricted $M$-torsors  is pointed with distinguished point
$(M, \Delta _M, \id _M)$.
Two  restricted $M$-torsors $(X, \Delta _X, \beta)$ and $(X', \Delta _{X'}, \beta')$ are {\sl equivalent}
if there  exists  $f \in \Aut _S(M)$ such that the composition {\hbox{$\beta \circ f \circ \beta '^{-1}: X' \lr X$}} is a morphism
of $(H,S)$-Hopf modules.
Denote by $\Torsr (M)$ the set of equivalence classes of  restricted {\hbox{$M$-torsors}}; it is pointed with distinguished point
the class of $(M, \Delta _M, \id _M)$. By Theorem 3.4 and Theorem 2.6 of the present article, by Proposition 2.8 and Theorem 3.1 of [8], one deduces the following statement:

\medskip

\noindent {\tencmb Corollary  3.9.} {\sl  Let $H$ be a Hopf algebra, $S$ be an $H$-comodule algebra, and $M$ be an $(H,S)$-Hopf module such that Condition $({\cal F}_n)$ is satisfied for  $0 \leq n \leq 2$.
Then there is a bijection of pointed sets
$$  \Torsr(M) \cong \Tors (H, \End_S(M)) .$$}

This result shows that, under weak technical conditions on $M$, the possible structures of $(H, \End_S(M))$-Hopf module on $\End_S(M)$ are closely related to the 
possible $(H,S)$-Hopf-module structures on $M$. More precisely, if $\End _S(M)$ is equipped with an $(H, \End_S(M))$-Hopf module structure $\Delta $, then following the track of
$\Delta $ along the four  isomorphisms $$\Tors (H, \End_S(M)) \cong \HC^1(H, \End_S(M)) \cong \Hr^1(H, M) \cong \Torsr(M), $$ one gets an $(H,S)$-Hopf-module structure $\Delta '$ on $M$ defined on 
an element $m \in M$ by
$$\Delta '(m) = \varphi _0(m_0) \ot \varphi_1 m_1.$$
Here we denote by $\varphi _0 \ot \varphi_1$ the element $\Delta (\id _M) \in \End_S(M)\ot H$, and as usual, we adopt the convention $\Delta _M(m) = m_0\ot m_1$.

\bigskip
\bigskip
\bigskip
{
\noindent {R{\eightrom EFERENCES}}

\medskip
\medskip

\item{[1]} A. B{\eightrom LANCO} F{\eightrom ERRO}, 
Hopf algebras and Galois descent, 
{\it Publ. Sec. Mat. Universitat  Aut\`onoma Barcelona} {\tencmb 30} ({\it 1986}), n$^{\hbox{{\fiverom o}}}$ $\!$1, 65 -- 80.
\medskip

\item{[2]} Ch. K{\eightrom ASSEL},  
{\it Quantum Groups},
 Graduate  Texts in Mathematics 155, 
Springer-Verlag,   New York ({\it1995}). 
\medskip

\item{[3]} H. F. K{\eightrom REIMER}, M. T{\eightrom AKEUCHI}, Hopf algebras and 
Galois extensions of an algebra,
{\it Indiana Univ. Math. J.} 
{\tencmb 30} ({\it1981}), n$^{\hbox{{\fiverom o}}}$ $\!$5,
675  --  692.
\medskip

\item{[4]}  S. L{\eightrom ANG}, J. T{\eightrom ATE}, Principal homogeneous spaces over abelian varieties,
{\it Amer. J. Maths.}   {\tencmb 80} ({\it1958}), 
{\hbox{659  --  684.}}
\medskip

\item{[5]} L. {\eightrom LE} B{\eightrom RUYN}, M. {\eightrom VAN DEN}  B{\eightrom ERGH}, 
F. {\eightrom VAN} O{\eightrom YSTAEYEN},
{\it Graded orders}, Birkh\"auser, Boston -- Basel ({\it1988}).
\medskip

\item{[6]} J.-L. L{\eightrom ODAY},
{\it Cyclic homology},  Grundlehren der Mathematischen Wissenschaften 301, Springer-Verlag, Berlin ({\it1988}).
\medskip

\item{[7]} Ph. N{\eightrom USS},   Noncommutative descent and non-abelian cohomology, 
{\it{\it K}-Theory} {\tencmb 12} ({\it1997}), n$^{\hbox{{\fiverom o}}}$ $\!$1, 23 -- 74.
\medskip

\item{[8]} Ph. N{\eightrom USS}, M.  W{\eightrom AMBST}, 
Non-Abelian Hopf Cohomology, {\it J. Algebra}
{\tencmb 312}, ({\it2007}), n$^{\hbox{{\fiverom o}}}$ $\!$2, 733 -- 754.
\medskip

\item{[9]} J.-P. S{\eightrom ERRE}, 
{\it Corps locaux}, Troisi\`eme \'edition corrig\'ee,
Hermann, Paris  ({\it1968}).

\medskip

\item{[10]} J.-P. S{\eightrom ERRE}, 
{\it Galois cohomology}, Springer-Verlag, Berlin --
Heidelberg  ({\it1997}). Translated from
{\it Cohomologie galoisienne},
Lecture Notes in Mathematics 5, Springer-Verlag, Berlin --
Heidelberg -- New York ({\it1973}).

\medskip

\item{[11]} M. E. S{\eightrom WEEDLER},  
 Cohomology of algebras over Hopf algebras, {\it Trans. Amer. Math. Soc.} {\tencmb 133}  ({\it1968}), 205~--~239.

}

\vfill
\vfill

\bye